\documentclass[10pt,reqno]{amsart}

\usepackage{amsmath, amsfonts, amssymb, amscd, enumerate}
\theoremstyle{plain}
\newtheorem{theo}{Theorem}[section]
\newtheorem{lem}[theo]{Lemma}
\newtheorem{prop}[theo]{Proposition}

\newtheorem{cor}[theo]{Corollary}
\theoremstyle{definition}
\newtheorem{rem}[theo]{Remark}

\newtheorem{definition}[theo]{Definition}

\newenvironment{pf}{\noindent{\it Proof. }}{$\square$\par\medskip}

\theoremstyle{plain}
\newtheorem{lemma}[theo]{Lemma}

\theoremstyle{definition}

\renewcommand{\=}{:=}

\newcommand{\beq}{\begin{equation}}
\newcommand{\eeq}{\end{equation}}
\renewcommand{\a}{\alpha}
\renewcommand{\b}{\beta}

\newcommand{\g}{\gamma}

\renewcommand{\l}{\lambda}
\newcommand{\m}{\mu}
\renewcommand{\o}{\omega}

\renewcommand{\r}{\rho}
\newcommand{\s}{\sigma}
\renewcommand{\t}{\tau}

\renewcommand{\L}{\Lambda}

\newcommand{\B}{\mathcal B}

\newcommand{\bC}{\mathbb{C}}
\newcommand{\bR}{\mathbb{R}}

\newcommand{\bH}{\mathbb{H}}

\newcommand{\bO}{\mathbb{O}}

\renewcommand{\gg}{\mathfrak{g}}
\newcommand{\gh}{\mathfrak{h}}
\newcommand{\gk}{\mathfrak{k}}
\newcommand{\gl}{\mathfrak{l}}
\newcommand{\gm}{\mathfrak{m}}
\newcommand{\gn}{\mathfrak{n}}

\newcommand{\gp}{\mathfrak{p}}

\newcommand{\gs}{\mathfrak{s}}
\newcommand{\gt}{\mathfrak{t}}

\newcommand{\su}{\mathfrak{su}}

\newcommand\GL{\mathrm{GL}}
\newcommand\SL{\mathrm{SL}}
\newcommand\SO{\mathrm{SO}}
\newcommand\SU{\mathrm{SU}}
\newcommand\U{\mathrm{U}}
\newcommand\T{\mathrm{T}}

\newcommand\Sp{\mathrm{Sp}}
\renewcommand\sp{\mathfrak{sp}}

\newcommand\G{\mathrm{G}}

\newcommand{\cB}{\mathcal{B}}
\newcommand{\cC}{\mathcal{C}}

\newcommand{\cI}{\mathcal{I}}

\newcommand{\cL}{\mathcal{L}}
\newcommand{\cM}{\mathcal{M}}

\newcommand{\cP}{\mathcal{P}}
\newcommand{\cQ}{\mathcal{Q}}

\newcommand{\cT}{\mathcal{T}}
\newcommand{\cU}{\mathcal{U}}
\newcommand{\cV}{\mathcal{V}}

\renewcommand{\square}{\kern1pt\vbox
{\hrule height 0.6pt\hbox{\vrule width 0.6pt\hskip 3pt
\vbox{\vskip 6pt}\hskip 3pt\vrule width 0.6pt}\hrule height0.6pt}\kern1pt}

\DeclareMathOperator\End{End\;}

\DeclareMathOperator\Ad{Ad}
\DeclareMathOperator\ad{ad}

\DeclareMathOperator{\Span}{Span}

\renewcommand\Im{\operatorname{Im}}

\newcommand{\Hom}{{\operatorname{Hom}}}

\newcommand{\wt}{\widetilde}
\newcommand{\wh}{\widehat}

\newcommand{\n}{\nabla}

\renewcommand{\mod}{\operatorname{mod}}
\newcommand{\be}{\begin{equation}}
\newcommand{\ee}{\end{equation}}

\def\<#1,#2>{\langle\,#1,\,#2\,\rangle}
\newcommand{\arr}{\begin{array}{rlll}}
\newcommand{\ea}{\end{array}}
\newcommand{\bea}{\begin{eqnarray}}
\newcommand{\eea}{\end{eqnarray}}
\newcommand{\bean}{\begin{eqnarray*}}
\newcommand{\eean}{\end{eqnarray*}}

\newcommand{\fff}{{\mathfrak f}_4}

\def\sideremark#1{\ifvmode\leavevmode\fi\vadjust{
\vbox to0pt{\hbox to 0pt{\hskip\hsize\hskip1em
\vbox{\hsize3cm\tiny\raggedright\pretolerance10000
\noindent #1\hfill}\hss}\vbox to8pt{\vfil}\vss}}}

\newcounter{ssig}
\newcounter{ttig}
\begin{document}
%
%
\title[Six-dimensional nearly K\"ahler manifolds of cohomogeneity one (II)]{Six-dimensional nearly K\"ahler manifolds\\ of cohomogeneity one (II)}
\author{Fabio Podest\`a and Andrea Spiro}
\begin{abstract} Let  $M$ be a six dimensional manifold, endowed with a cohomogeneity one action of $G= \SU_2\times\SU_2$, and   $M_{\text{reg}} \subset M$ its subset of regular points.  We show that   $M_{\text{reg}}$  admits  a
smooth, 2-parameter family of $G$-invariant,  non-isometric strict nearly K\"ahler  structures and   that a  1-parameter subfamily of such structures smoothly extend over a singular orbit of type $S^3$. This determines  a new class  of examples of nearly K\"ahler structures on $TS^3$.  \end{abstract}
\maketitle

\section{Introduction}
\setcounter{equation}{0}
\bigskip

A Riemannian manifold $(M,g)$ is called {\it nearly K\"ahler} (shortly NK) if it admits  an almost complex structure $J$,  such that  $g$ is Hermitian and the Levi-Civita connection $\n$ satisfies  $(\nabla_XJ)X = 0$ for every vector field $X$.  A NK structure $(g,J)$ is called {\it strict} if $\nabla_v J|_p\neq 0$ for every $p\in M$ and every $0\neq v\in T_pM$ (see  e.g. \cite{Gra1,Gra2, M, Na1, Na2} for main properties).
Significant examples  are  the so-called  $3$-symmetric spaces  with their canonical almost complex structures. Recall also that the NK manifolds constitute one of the  sixteen classes of Gray and Hervella' s classification of almost Hermitian manifolds and their canonical Hermitian connection $D$ has totally skew and $D$-parallel torsion (\cite{Gra3,GH, Ag,Na1,FI}).\par
According to Nagy's structure theorem (\cite{Na1,Na2}), any complete strict NK manifold
is finitely covered by a product of homogeneous $3$-symmetric manifolds, twistor spaces of positive quaternion
K\"ahler manifolds with their canonical NK structure and six dimensional strict NK manifolds.  This is one of the reasons which raise a particular  interest for   six dimensional strict NK structures. \par
It is also known  that, in six dimensions, the ``strictness''  condition is equivalent to the fact that the NK
structure is not K\"ahler and that strict  NK manifolds are automatically Einstein and related  with the existence of  a nonzero
Killing spinor (see e.g. \cite{Gra2, Gru}). Other reasons of interest for   NK structures in six dimensions are provided  by their relations with
geometries with torsion,  $\G_2$-holonomy and supersymmetric models (see e.g. \cite{Ag,Na3, FI, Ba}). \par
\smallskip
Up to now, the only known examples of compact, six dimensional, strict NK manifolds are the six dimensional $3$-symmetric spaces
endowed with their natural NK structures, namely the standard sphere $S^6 = \G_2/\SU_3$, the twistor spaces
$\bC P^2 = \Sp_2/\U(1)\times \Sp_1$ and $F = \SU_3/\U(1)^2$, and the space $S^3 \times S^3 = \SU_2^3/\SU_2$. This is    actually   the  list of   homogeneous strict NK manifolds in six dimensions (\cite{Bu}). \par
\par
In \cite{PS},  we started the study of six dimensional strict NK manifolds $(M, g, J)$, admitting a compact
Lie group $G$ of automorphisms with a codimension one orbit. In these hypothesis, we proved  that if $M$ is compact:
\begin{itemize}
\item[a)]  $G$ is  semisimple and locally isomorphic to $\SU_3$ or $\SU_2 \times \SU_2$;
\item[b)] if $G$ is locally isomorphic to $\SU_3$, then $(M,g)$ has constant sectional curvature (this is true also if $M$ not compact);
\item [c)] if $G$ is locally isomorphic to  $\SU_2 \times \SU_2$, then
$M$ is $G$-diffeomorphic to one of the following $3$-symmetric spaces: the sphere $S^6 =\G_2/\SU_3$, the
 projective space $\bC P^2 = \Sp_2/\U(1)\times \Sp_1$ or $S^3 \times S^3 = \SU_2^3/\SU_2$.
\end{itemize}
\par
In this paper we focus on the six dimensional
strict NK manifolds $(M,g,J)$ (not necessarily compact),   on which   $G = \SU_2\times\SU_2$ acts by automorphisms with cohomogeneity one, i.e. with a codimension one orbit. In particular, we know that a principal isotropy subgroup is isomorphic to a one-dimensional torus $T^1_{\operatorname{diag}}$, which is diagonally embedded in $G$ (see \cite{PS}). \par
We  start from  the following known fact: any strict NK structure $(g, J)$ on a given six dimensional,  oriented manifold   can be completely recovered  by its K\"ahler form $\o = g(J \cdot, \cdot)$ and any non-degenerate $2$-form $\o$, whose differential $d\omega$ is stable in the sense of Hitchin (\cite{Hi1,Hi2}) and satisfying  a suitable differential problem, is the K\"ahler form of a strict NK structure $(g, J)$ (see  Theorem \ref{startingpoint} for the complete  statement).  From this,  we  obtain that,  locally, any  $G$-invariant NK structure   is completely determined by  a  smooth map $f = (f_1, \dots, f_5): ]a, b[ \subset \bR \longrightarrow \bR^5$, which solves a certain differential problem and which
 gives  the components of the K\"ahler form  $\o$ with respect to a special basis of $G$-invariant $2$-forms along the points of a normal geodesic $\g: ]a,b[ \longrightarrow M$.  We then study  the solvability of the differential problem on  $f$, first on an open set of $G$-principal points and then on a suitable tubular neighborhood of a singular orbit diffeomorphic to $S^3$. Our main results can be   outlined as follows.\par

\begin{theo} \label{mainresults} Let $G = \SU_2 \times \SU_2$. Then:\par
\begin{itemize}
\item[1)] There exists a  2-parameter family  of non-isometric,  non locally homogeneous, $G$-invariant strict  NK structures on $G/K \times \bR$, with $K = T^1_{\operatorname{diag}}$;
\item[2)]  There exists a   family of non isometric,  $G$-invariant strict NK structures on $TS^3 \cong G \times_{\SU_2} \bR^3$,  smoothly parameterized  by   $]0, + \infty[ \subset \bR$.
All such structures  are non locally homogeneous, except precisely two of them, which  are $G$-equivalent to  those of suitable tubular neighborhoods of the singular $G$-orbit  $S \simeq S^3$  in  $S^6$ and $S^3 \times S^3$, respectively.
\end{itemize}
\end{theo}
As an interesting consequence of  (2), we have that the two locally homogeneous NK structure $(g_0, J_0)$ and $(g_1,J_1)$ on $TS^3$ are connected by a smooth family of non locally homogeneous NK structures $(g_t,J_t)$, $t \in ]0,1[$. The interesting problem to determine how many of these NK structures can be $G$-equivariantly completed remains unsolved. Note that any such structure admitting a $G$-equivariant compactification would give a new NK structure on either
$S^6$ or $S^3\times S^3$.  \par
\smallskip
Other interesting information are given by  the proofs of the above results. In fact,  the  proof of (1) clearly  shows that any smooth  family of non-isometric $G$-invariant
NK structures on  $G/K \times \bR$ can be described  by at most two parameters. The proof of (2) shows that  any $G$-invariant NK structure $(g, J)$ on $T S^3$ is isometric to one of the structures described above. We may therefore say  that  the isometric  moduli space $\cM$ of such  NK structures on $T S^3$ can be identified with  $\cM = \bR$.\par
The structure of the paper  is as follows. In \S 2 and \S 3, we recall some known facts on stable 3-forms and six dimensional  strict NK structures and use them to determine the differential problem that characterizes the map $f = (f_1, \dots, f_5)$,  associated with  $G$-invariant strict NK structures. In \S 4, we show that, modulo a finite group of transformations,  the map  $f $ is an isometric invariant of  any such  NK structure and we give the explicit expressions for  the $f$'s associated with  locally homogeneous NK structures. In \S 5 and \S 6 we prove claim (1) and (2) of Theorem \ref{mainresults}, respectively. \par
\smallskip
Given a semisimple Lie group $G$ (which, in almost all the paper, is always assumed to be  $G = \SU_2 \times \SU_2$)  and an action of $G$ on a manifold $M$, we denote by $\gg$ its Lie algebra,   by $\cB$  the  Cartan-Killing form of $\gg$ and we use symbols of the form    $\wh X$  to indicate the vector fields on $M$ corresponding to elements   $X \in \gg$.\par
\bigskip
\bigskip
\bigskip
\section{Preliminaries}
\setcounter{equation}{0}
\bigskip
\subsection{Stable 3-forms  of oriented  $6$-dimensional vector spaces} \label{preliminary1}
Let  $V$ be a $6$-dimensional real vector space  with a fixed orientation and consider the standard action on  $\Lambda^3V$ of the orientation preserving transformation in  $G = \GL^+(V)$.  It is known (e.g. \cite{SK, Hi2}) that  $\Lambda^3V$ can be divided  into the following disjoint $G$-invariant  sets: a  $G$-invariant hypersurface $(\L^3 V)_{0} = \{\ P = 0\ \}$, given by  the zero set of
 a suitable relative $G$-invariant  irreducible polynomial  $P$ of degree $4$,  and the  complementary   open  sets
$$(\L^3 V)_{-} = \{\ P < 0\ \}\ ,\qquad (\L^3 V)_{+} = \{\ P > 0\ \} \ .$$
These two sets are  open $G$-orbits and their generic stabilizers have connected components conjugate to $\SL_3(\bC)$ and $\SL_3(\bR)\times \SL_3(\bR)$ respectively.\par
The polynomial $P$ is defined as follows.  We consider the isomorphism $A: \Lambda^5 V^*\longrightarrow  \Hom(V^*,\Lambda^6 V^*) \cong V\otimes \Lambda^6 V^*$ that is induced by the wedge product $\wedge: \Lambda^5V^* \otimes V^*\to \Lambda^6V^* $. We fix a non zero element $\tau\in \Lambda^6V^*$ and for every $\theta\in \Lambda^3V^*$ we define $S_\theta\in \End(V)$ as follows: given $v\in V$

\beq \label{2.1} A(\imath_v \theta \wedge \theta) = S_\theta(v)\otimes \tau.\eeq
One can check that  $S_\theta^2 = P(\theta) Id_V$  for some polynomial map $P: \Lambda^3 V^* \to \bR$ and
that $P$ is irreducible, of degree $4$ and $\SL(V)$-invariant (see e.g. \cite{SK}, p.80), which depends  on  the choice of volume form as follows: if $P$, $P'$ are determined using volume forms $\t$,  $\t' = c \cdot \t$, respectively, then
\beq \label{star} P' = \frac{1}{c^2} P \ . \eeq
\bigskip
By definitions, for any  $\theta\in (\L^3 V)_{-} = \{\ P < 0\ \}$, the endomorphism
\begin{equation}\label{J} J_\theta \= \frac{1}{\sqrt{-P(\theta)}} S_\theta\end{equation}
is a complex structure.   We call any  3-form $\theta\in (\L^3 V)_{-}$ a {\it stable 3-form\/} and the corresponding
$J_\theta$ the {\it complex structure determined by  $\theta$ (relatively to the  given orientation)}. We conclude pointing out that for any $\theta \in (\L^3 V)_{-} $ the  complex 3-form
\beq  \a= \frac{1}{2}\left(\theta + i J^*_{\theta} \theta\right)\eeq
is  of type (3,0) (see e.g.  Prop. 2 and formulae (8), (9) in \cite{Hi2}). In particular,  $\imath_{(v + i J_\theta v)} \a = 0$ for any $v \in V$ and hence, for any $v_1, v_2, v_3 \in V$,
\beq \label{furbata}Ê\theta(J_\theta v_1, J_\theta v_2, J_\theta v_3) =  \theta(J_\theta v_1, v_2, v_3)\ .\eeq
\subsection{Stable 3-forms and NK-structures} Let $M$ be a $6$-dimensional oriented manifold.
\begin{definition}  A 3-form $\psi$ on $M$ is called {\it stable\/} if $\psi_x $ is stable for any $x \in M$.
If $\psi$ is stable, we consider the almost complex structure $J_\psi$ on $M$ such that
$(J_\psi)|_x := J_{\psi_x}$ at every $x\in M$, where $J_{\psi_x}$ is defined as in (\ref{J}).\end{definition}
The following is well-known (see e.g. \cite{Re, Hi1, Hi2, CS, Sa}).\par
\begin{theo} \label{startingpoint} Let $\o \in \L^2 T^*M$ and $\psi \in \L^3 T^* M$ so that:
\begin{itemize}
\item[i)] $\psi$ is stable, $\o$ is $J_\psi$-invariant and  $g = \o(\cdot, J_\psi \cdot)$ is positively or negatively defined;
\item[ii)] there exists $ \mu\in \bR^+$ so that
\beq\label{fundamentalequations} \left\{\begin{array}{l}  d \o  = 3\psi \\
\phantom{a} \phantom{a}\\ d (J^*_\psi \psi)  =  - 2 \mu \cdot \o \wedge \o\ .\end{array} \right.\eeq
\end{itemize}
Then $(g, J_\psi)$ (or $(g, -J_\psi)$) is a strict NK structure on $M$ with  scalar curvature $s = 30 \mu$.\par
Conversely, let  $(g, J)$ be a strict NK structure on $M$ with (constant) scalar curvature  $s$
and denote by $\o = g(J \cdot, \cdot)$ and $\psi = \frac{1}{3} d \o$.  Then  $\o$ and $\psi$ satisfy (i)
and (ii) with $ \m = \frac{s}{30}$ and $J = \pm J_\psi$.
\end{theo}
\begin{rem} In (i), the condition ``$\o$ is $J_\psi$-invariant'' is indeed redundant, if the system (\ref{fundamentalequations})
is satisfied. Namely, $3\omega\wedge \psi = \frac{1}{2}d(\omega^2) = 0$.
\end{rem}
\begin{rem} Conditions (i) and (ii) were considered for the first time by Reyes Carri\'on in his PhD thesis.
\end{rem}
\bigskip
\subsection{NK-structures of cohomogeneity one for an $\SU_2 \times \SU_2$-action} \label{notation}  We recall that the
action of  a compact connected Lie group $G$, acting almost effectively and isometrically on a Riemannian
manifold $(M,g)$, is of {\it cohomogeneity one} if the generic orbits have codimension one.  The  points in such
generic orbits (i.e. whose $G$-isotropy is, up to conjugation,  minimal)
 are called {\it regular} and constitute an open and dense  $G$-invariant subset $M_{\operatorname{reg}}  \subset M$.
 \par
If we  denote by $\xi$ a unit vector field on $M_{\operatorname{reg}}$, which is orthogonal to all $G$-orbits, its integral curves are geodesics  that meet every $G$-orbit orthogonally. Moreover, any regular orbit $G\cdot p = G/K$ admits a tubular neighborhood which is $G$-equivariantly isometric to
the product $]a,b[\times G/K$, $]a,b[\subset \bR$, endowed with the metric
\beq \label{2.7}  g = dt^2 + g_t,\eeq
where $g_t$ is a smooth family of $G$-invariant metrics on $G/K$ and the vector field $\partial/\partial t$ corresponds
to $\xi$. \par
Let $(g, J)$ be an NK structure on $M$ and $G$ a compact connected Lie group acting isometrically  with cohomogeneity one
on $(M, g)$.
Throughout the following we will always suppose that $G$ {\it preserves the almost complex structure $J$\/}. Note that this condition is automatically
 satisfied if the manifold is compact and the metric has not constant sectional
curvature (see \cite{M}, Prop. 3.1).  \par
Due to \eqref{2.7},  any  regular point of an NK manifold  of cohomogeneity one  admits a neighborhood which is
locally identifiable with a Riemannian manifold of  the form
 $M = ]a,b[\times G/K$,  endowed with the NK-structure $(g, J)$ associated to
a $G$-invariant  pair $(\o, \psi)$, that satisfies (i) and (ii) of  Theorem \eqref{startingpoint}
 and  so that for every $X\in \gg$
 \beq \label{2.10} g\left(\xi, \xi\right)  = 1\ ,\ \ \
 g(\xi,\hat X)  = 0.\eeq
In the rest of the paper, we will be concerned with the strict NK-structures on $6$-dimensional manifolds of
cohomogeneity one w.r.t. $G = \SU_2 \times \SU_2$.  By Lemma 3.3 in \cite{PS}, which holds also when the manifold is not compact, we can suppose that, up to an automorphism of $\gg$, the Lie algebra $\gk$ of a regular isotropy subgroup $K$ is diagonally embedded into a Cartan
subalgebra of $\gg$, i.e.
\beq\label{2.11} \gk =  \bR\cdot (H, H)\in \su_2 + \su_2\eeq
for some generator $H$ of a Cartan subalgebra $\gh \subset \su_2$. \par
\medskip
We conclude this section by fixing some notation and a particular  basis for  $\su_2 + \su_2$,  which  is particularly useful for  our computations and will be constantly used in the rest of the paper.  \par
\begin{itemize}
\item[--] $\cB$ always denotes the (negatively defined) Cartan-Killing form of $\su_2$;
\item[--] $\gn = (\bR\cdot H)^\perp$ is  the $\cB$-orthogonal
complement of $\gh = \bR\cdot H$ in $\su_2$;
\item[--] with no loss of generality, we always assume that  $\ad(H)|_\gn$ is  a complex structure on $\gn$;
\item[--] $(E, V)$ is a basis for $\gn$  with  $\cB(E, E) = \cB(V,V) = -1$ and $V = [H, E]$;
\item[--]  the basis of $\su_2 + \su_2$ which we constantly consider is  given by the elements
$$ U \= (H, H)\ ,\quad A \= (H, - H)\ ,$$
\beq \label{basissu2} E_1 \= (E,0)\ ,\quad V_1 \= (V, 0)\ ,\quad E_2 \= (0, E)\ ,\quad V_2 \= (0,V)\ ;\eeq
\item[--]  $\g$ is the curve $\g_t = (t, eK) \in ]a,b[ \times (\SU_2 \times \SU_2)/K$ and
$\xi = \frac{\partial}{\partial t}$;
\item[--] for any point $\g_t$, we denote by $B_t$  the basis for  $ T_{\g_t} M$ equal to
\beq\label{basisBt} B_t =  (\xi, \wh A, \wh E_1, \wh V_1,  \wh E_2, \wh V_2)_{\g_t}  \ ; \eeq
we also denote  by  $B^*_t = (\xi^*, A^*, E_1^*, V^*_1, E^*_2, V^*_2)_{\g_t}$  the corresponding dual
coframe in $T^*_{\g_t} M$.  With no loss of generality, we will always assume that $M$ is oriented and
that $B_t$ is in the prescribed orientation  for any $t$.
\end{itemize}

\section{The equations}
\setcounter{equation}{0}
\label{theequations}
\medskip
In this section, we want to determine the differential problem that  characterizes the $G$-invariant
pairs $(\o, \psi)$, with $G = \SU_2 \times \SU_2$,  on a manifold $M$ as in \eqref{2.11}, corresponding
to NK structures. \par
\medskip
We first describe the  $G$-invariant 2- and 3-forms on $M$. With the notation and assumptions
of  \S  \ref{notation}, any  $G$-invariant p-form $\varpi$ is uniquely determined by
the values of $\varpi_{\g_t}$ on the tangent spaces $T_{\g_t} M$, $\g_t = (t, eK)$.
Since the curve $\g(]a,b[) \subseteq Fix(K)$, the form $\varpi$ is $G$-invariant if and only
if $\varpi_{\g_t}$ is $K$-invariant for any $t$. \par
Now,  the tangent space  $T_{\g_t} M$ decomposes  into the following $K$-moduli
\beq \label{decomposition} T_{\g_t} M = < \wh E_1|_{\g_t},\wh V_1|_{\g_t}> \oplus <\wh E_2|_{\g_t},\wh V_2|_{\g_t}>
\oplus <\xi_{\g_t}, \hat A_{\g_t} > \ ,\eeq
where $K$ acts irreducibly on the first  two and trivially on the last one. \par
 Using  this, we determine the spaces of $K$-invariant 2-forms of the tangent spaces $T_{\g_t}M$ and we obtain  that   the  space of $G$-invariant $2$-forms is generated over $C^{\infty}(]a,b[)$ by
the five invariant  forms $\o^i$, $1 \leq i \leq 5$,  defined by
\par
\beq \label{invariantforms}
\omega^1|_{\g_t} = \xi^*\wedge A^*\ ,\qquad \omega^2|_{\g_t} = E_1^*\wedge V_1^*\ ,\qquad\ \omega^3|_{\g_t} = E_2^*\wedge V_2^*\ ,$$
\par
$$\omega^4|_{\g_t} = E_1^*\wedge E_2^* + V_1^*\wedge V_2^*\ ,\qquad \omega^5|_{\g_t} = E_1^*\wedge V_2^*- V_1^*\wedge E_2^*\ .\eeq
\par
\noindent
In a similar way,  we obtain that  the space of $G$-invariant $3$-forms on $M$ is generated over $C^{\infty}(]a,b[)$ by
the eight invariant  $3$-forms $\psi^{a i}$, $a = 1,2$,  $i = 2, \dots, 5$,  defined by
$$\psi^{1i} |_{\g_t} \=  \xi^*\wedge \omega^i|_{\g_t} \ ,\qquad \psi^{2 i}|_{\g_t}  \= A^*\wedge \omega^i|_{\g_t}\
\quad 2 \leq i \leq 5\ .$$
Due to this,  the $G$-invariant strict NK structures on $M$ are in natural correspondence  with  the collections of real functions $f_i, p_{aj}Ê\in C^{\infty}(]a,b[)$ such that the  pairs
\begin{equation}\label{opsi} \left(\omega  =    f_i \omega^i\ \  ,\ \  \psi  =  p_{aj}  \psi^{aj}\right)\end{equation}
satisfy  the conditions of Theorem \ref{startingpoint} together with the constraints \eqref{2.10}. \par
\medskip
 Let us now consider the  following lemma.\par
\begin{lem} \label{lemma1} The condition $\psi = \frac{1}{3} d \o$ on forms \eqref{opsi} is equivalent to the equations
\beq \label{3.2}
 p_{22} = p_{23}Ê= 0\ ,\eeq
\beq  \label{3.3} p_{12} = \frac{f_2'}{3} + \frac{f_1}{12} \ ,\qquad p_{13} =  \frac{f_3'}{3}- \frac{f_1}{12}\ ,\eeq
\beq  \label{3.4} p_{14} = \frac{f'_4}{3} \ ,\qquad p_{15} = \frac{f'_5}{3} \ ,\qquad p_{24} = \frac{2}{3} f_5\ ,\qquad
p_{25} = - \frac{2}{3} f_4\ .\eeq
\end{lem}
\begin{pf} Using the fact that  the flow $\Phi^\xi$ of $\xi = \frac{\partial}{\partial t}$ commutes with the action
of $G$ on $M$, it is immediate to realize that the  $G$-invariant  forms $\o^i$ are also $\Phi^\xi$-invariant. By
$G$- and $\Phi^\xi$-invariance  and Koszul's formula, it follows  that for any $X, Y, Z \in \su_2 + \su_2$
\beq  \label{differential1}Ê d \o^i(\xi, \wh X, \wh Y) =  \o^i(\xi, \wh{[X,Y]})\eeq
and
\beq \label{differential2} d \o^i(\wh X, \wh Y, \wh Z) = \o^i( \wh X, \wh{[Y,  Z]}) +
\o^i( \wh Y, \wh{[Z, X ]}) +  \o^i( \wh Z, \wh{[X, Y ]}) \ .\eeq
Using \eqref{differential1}, \eqref{differential2} and the fact that
$[E_i, V_i] = \frac{(-1)^{i+1}}{4}A$ $(\mod \bR U)$ for $i=1,2$,
we see that
\beq \label{differentials1} d \o^1_{\g_t} = \frac{1}{4}\xi^* \wedge  \left(\o^2 - \o^3\right)|_{\g_t} \ , \quad
d\o^2_{\g_t} = d\o^3_{\g_t} = 0\ ,\eeq
\beq \label{differentials2}Ê d\o^4_{\g_t} =  -2 (A^*\wedge \omega^5)|_{\g_t},\ \qquad
d\omega^5_{\g_t} = 2 (A^*\wedge \omega^4)|_{\g_t}\ .\eeq
From this and the $G$-equivariance, we have that the equality
$$  \psi = p_{1 i} \xi^* \wedge \o^i + p_{2i}ÊA^* \wedge \o^i =   \frac{1}{3} d \o = \frac{1}{3} d(f_j\o^j)$$
 is equivalent to \eqref{3.2} - \eqref{3.3}. \end{pf}
The next lemma gives the conditions corresponding to the stability of $\psi$ and to condition \eqref{2.10}. \par
\begin{lem} \label{lemma2} Given a pair \eqref{opsi} satisfying \eqref{3.2} - \eqref{3.4}, we have that $\psi$ is stable,
$\o$ is $J_\psi$-invariant and $g = \o(\cdot , J_\psi \cdot)$ satisfies  \eqref{2.10}  if and only if at all
points of  $]a,b[$ the following conditions hold
\begin{itemize}
\item[i)] $f_1 < 0$  and $f_4$ and $f_5$ are of the form
\beq \label{firstconditionbis} f_4 =  \fff \cos \theta_o \ ,\qquad f_5 = \fff \sin \theta_o\ ,\eeq
for a suitable function  $0 < \fff \in \cC^\infty(]a,b[)$ and some constant $\theta_o \in \bR$;
\item[ii)]
\begin{equation} \label{firstcondition}
4 \fff^2 - \left((\fff')^2 - \left(f'_2 + \frac{f_1}{4}\right) \left(f'_3 - \frac{f_1}{4}\right)\right)
(f_1)^2 = 0  \ ,\end{equation}
\item[iii)]
\beq\label{stability}
(\fff')^2  - \left(f'_2 + \frac{f_1}{4}\right) \left(f'_3 - \frac{f_1}{4}\right) > 0\ .
\eeq
\end{itemize}
If (i) -(iii) are satisfied, $J_\psi$ is represented in the basis \eqref{basisBt} by a  matrix of the form
\beq \label{complexstructure}
J_\psi =  \left( \begin{matrix} K& 0 \\
0 & L
\end{matrix} \right)\eeq
where, if we put $q = p_{14}^2+p_{15}^2 - p_{12} p_{13}$,
\beq \label{K} K =  \left(\smallmatrix 0& \sqrt{\frac{p_{24}^2+p_{25}^2}{q}}\\
- \sqrt{\frac{q}{p_{24}^2+p_{25}^2}}
 &0
\endsmallmatrix\right)\ ,
\eeq
$$L=  \smallmatrix\phantom{a}\\ \frac{1}{\sqrt{q(p_{24}^2 + p_{25}^2)}}\endsmallmatrix\cdot
\left(\smallmatrix
0 & -  p_{15} p_{24} +  p_{14} p_{25}   &  p_{13} p_{25}& - p_{13} p_{24}\\
 p_{15} p_{24} -  p_{14} p_{25}  & 0 & p_{13} p_{24} &  p_{13} p_{25}\\
 p_{12}  p_{25} & p_{12}  p_{24} & 0 & p_{15} p_{24} -  p_{14} p_{25} \\
- p_{12} p_{24} & p_{12}  p_{25}  & -p_{15}  p_{24} +  p_{14} p_{25} & 0
\endsmallmatrix\right).$$

\end{lem}
\begin{pf} We fix a point $p = \g_{t_o}$ in the curve $\g$ and we consider  the orientation of $T_p M$ given by $B_{t_o}$.
From  definitions, up to a factor,  we get
\beq \label{3.15} S_{\psi_p} \left(\xi_p\right) = (p_{14} p_{24} +p_{15} p_{25})\xi_p + (p_{12} p_{13} -  p_{14}^2 - p_{15}^2) \hat A_p\ ,\eeq
\beq \label{3.16}  S_{\psi_p} (\hat A_p) = (p_{24}^2 + p_{25}^2)\xi_p - (p_{14} p_{24}  + p_{15} p_{25})\hat A_p\ .\eeq

Assume now that $\psi_p$ is stable (i.e. $P(\psi_p) < 0$) and that $\eqref{2.10}$ holds. Since $S_{\psi_p}= \sqrt{-P(\psi_p)} J_{\psi_p} $ and $g_p(J_{\psi_p}(\xi_p),\xi_p) = 0$, we have
\beq  \label{3.17} g(S_{\psi_p}(\xi_p),\xi_p) =  p_{14} p_{24} + p_{15} p_{25} = 0 \quad\overset{\eqref{3.4}} \Rightarrow\quad f'_{4} f_5  - f'_5  f_4 = 0\ .\eeq
Using \eqref{3.17} in \eqref{3.15} and \eqref{3.16}, we get that
$$P(\psi_p) = -(p_{24}^2+p_{25}^2)(p_{14}^2+p_{15}^2 - p_{12} p_{13}) = $$
\begin{equation}\label{pol} = - \frac{4}{81} \left(f_4^2 + f_5^2\right) \left((f'_4)^2 + (f'_5)^2 - \left(f'_2 + \frac{f_1}{4}\right) \left(f'_3 - \frac{f_1}{4}\right)\right) < 0\ .
\end{equation}
 From this we see that $f_4^2 + f_5^2 > 0$. By \eqref{3.17} either $\frac{f_4}{f_5}$ or $\frac{f_5}{f_4}$ is constant (according to the case $f_4 \neq 0$ or $f_5 \neq 0$). If we set $\fff \= \sqrt{f_4^2 + f_5^2}$, we get that $f_4$ and $f_5$ are as in \eqref{firstconditionbis}.
 Moreover, using \eqref{3.15}, \eqref{2.10} and \eqref{pol}
\beq\label{f1} f_1(t_o) =  \o_p(\xi, \wh A) = - \frac{\sqrt{-P(\psi_p)}}{ p_{14}^2 + p_{15}^2 - p_{12} p_{13}} <  0 \eeq
 and (i) follows.
 Then, (ii) follows from equality \eqref{f1} using (i) and \eqref{3.3} - \eqref{3.4}, while (iii) follows from
 (ii) using the fact that $f_1$ and $\fff$ do not vanish.\par
 Conversely,  it is immediate to check  that  \eqref{firstconditionbis} -  \eqref{stability}  imply that $\psi_p$ is stable and that $\eqref{2.10}$ is satisfied.\par
\smallskip
Finally, assume \eqref{firstconditionbis} -  \eqref{stability}.  Through a direct but lengthy calculation, one can check that the complex structure $J_\psi =  \frac{1}{\sqrt{-P(\psi)}} S_\psi$ is  of the form \eqref{complexstructure} at any point  $\g_t$. Using this, one can also check that $\o$ is  $J_\psi$-invariant.\end{pf}
\bigskip
\begin{theo} \label{theoequations} Let  $(\o,\psi)$ be as  in  \eqref{opsi}. It   satisfies  \eqref{2.10} and all conditions of Theorem \ref{startingpoint},  but the positivity of $g = \o(\cdot, J_\psi \cdot)$,
if and only if the $p_{aj}$'s are as in \eqref{3.2} - \eqref{3.4} and the functions $f_i$'s satisfy the following:
\begin{itemize}
\item[i)] $f_4$ and $f_5$ are of the form
$ f_4 =  \fff \cos \theta_o $ and $f_5 = \fff \sin \theta_o$
for some constant  $\theta_o \in \bR$ and a positive function $0 < \fff \in \cC^\infty(]a,b[)$;
\item[ii)] $f_1 < 0$ and  $
 (\fff')^2  - \left(f'_2 + \frac{f_1}{4}\right) \left(f'_3 - \frac{f_1}{4}\right) > 0$
  at all points;
\item[iii)]  $f_1$, $f_2$, $f_3$, $\fff$ satisfy the differential system
\beq
\label{finalsystem} \left\{\begin{array}{l}
[(f_2' + \frac{1}{4} f_1)f_1]'  + 12 \mu f_1f_2 = 0\ ,\cr
\ \cr
[(f_3' - \frac{1}{4} f_1)f_1]' + 12\mu f_1f_3 =  0\ , \cr
\ \cr
(\fff'f_1)'  -  4\frac{\fff}{f_1}  + 12\mu  f_1 \fff = 0\ , \cr
\ \cr
  f_1(f_2' - f_3' + \frac{1}{2}f_1) + 48 \mu (f_2f_3 -  \fff^2) =  0\ ,
\end{array}\right.\eeq
together with the algebraic condition
\beq \label{spippolo}
\left. 4 \fff^2 - \left((\fff')^2 - \left(f'_2 + \frac{f_1}{4}\right) \left(f'_3 - \frac{f_1}{4}\right)\right) (f_1)^2 \right|_{t=t_o} = 0\eeq
to be satisfied  at some $t_o \in ]a,b[$.
\end{itemize}
\end{theo}
\begin{pf}
After   Lemmas \ref{lemma1} and \ref{lemma2}, it remains to consider the equation $d \wh \psi = - 2 \mu \o \wedge \o$ with $\wh \psi = J^*_{\psi} \psi$. \par
By the $G$-invariance, $\wh \psi = J^*_\psi \psi$ is of the form $\hat\psi = \sum_{a j} \wh p_{aj} \psi^{aj}$
for suitable $\wh p_{aj}\in C^\infty(]a,b[)$. Using \eqref{furbata} and \eqref{complexstructure}, we see that for $j = 2, \dots, 5$,
\beq \wh p_{1j} =  \sqrt{\frac{p_{14}^2+p_{15}^2 - p_{12} p_{13}}{p_{24}^2+p_{25}^2}} p_{2j}\ ,\
\wh p_{2j} =  - \sqrt{\frac{p_{24}^2+p_{25}^2}{p_{14}^2+p_{15}^2 - p_{12} p_{13}}} p_{1j}\ .\eeq
On the other hand, if we assume that \eqref{2.10} is satisfied, then \eqref{f1} holds and therefore,
using \eqref{3.2} - \eqref{3.4},
\beq \label{3.2bis}
\wh p_{21} = \wh p_{11} = \wh p_{12} = p_{13}Ê= 0\ ,\eeq
\beq  \label{3.3bis} \wh p_{22} =  \frac{f_1}{3}\left(f_2' + \frac{f_1}{4}\right) \ ,\qquad \wh p_{23} =  \frac{f_1}{3} \left(f_3' - \frac{f_1}{4}\right)\ ,\eeq
\beq  \label{3.4bis} \wh p_{24} = \frac{f_1 f'_4}{3} \ ,\qquad \wh p_{25} = \frac{f_1 f'_5}{3} \ ,\qquad \wh p_{14} = - \frac{2}{3} \frac{f_5}{f_1}\ ,\qquad \wh p_{15} =  \frac{2}{3} \frac{f_4}{f_1}\ .\eeq
We now compute $d \wh \psi$ along the curve $\gamma$. It is not difficult to see that
\beq   \label{differentials3}   d A^* =  \frac{1}{4} (\o^3 - \o^2) \eeq
and therefore, using  \eqref{differentials1},  \eqref{differentials2}, \eqref{differentials3}  together with
 the fact that $\omega^i \wedge \omega^j = 0$ for $j>i = 2,3$,
we get that,  at the points of $\g$,
\beq\label{dhatpsi}  d\hat\psi  =  \wh p_{22}' \xi^* \wedge A^* \wedge \o^2 +  \wh p_{23}'   \xi^* \wedge A^* \wedge \o^3 + \left( \wh p_{24}'  -   2 \wh p_{15}\right) \xi^*\wedge A^*\wedge \omega^4 +$$
$$
+ \left( \wh p_{25}' + 2  \wh p_{14} \right) \xi^*\wedge A^*\wedge \omega^5 +  \frac{1}{4}\left(  \wh p_{22} - \wh p_{23} \right) \xi^* \wedge \o^2 \wedge \o^3. \eeq
Now
\beq \label{omega2}\omega \wedge \o|_{\g_t}   =  2 \sum_{i=2}^5 f_1f_i\ \xi^*\wedge A^*\wedge\omega^i|_{\g_t}  + 2 ( f_2f_3 - f_4^2 - f_5^2)\ \omega^2\wedge\omega^3|_{\g_t} \eeq
and therefore, comparing \eqref{dhatpsi} and \eqref{omega2} and using \eqref{3.2bis}, \eqref{3.2bis}, \eqref{3.4bis}, we see that the equation $d \wh \psi = - 2 \mu \o \wedge \o$ is equivalent to the system of equations \eqref{finalsystem}.\par
Summing up, by  Lemmas \ref{lemma1} and \ref{lemma2} and previous arguments, the conditions of Theorem \ref{startingpoint} (with the only exception of the positivity of $g$) are equivalent to (i), (ii), \eqref{finalsystem}  and  \eqref{firstcondition}.
On the other hand, a straightforward  check shows that whenever the $f_i$'s satisfy the first three equations of \eqref{finalsystem}, then the derivative of  \eqref{firstcondition} coincides with
the derivative  of the last equation of \eqref{finalsystem}  multiplied by  $- \frac{f_1(t)^2}{4}$. Hence, if we assume that $f_1(t) < 0$, for any solution of \eqref{finalsystem}  the derivative of  \eqref{firstcondition} is identically  $0$ and
 \eqref{firstcondition} is satisfied as soon as it is satisfied at just one point $t_o = ]a,b[$.\end{pf}
\begin{rem} \label{succinilcolina}Given a NK structure ($g,J$) on the manifold $M$, we may consider the
corresponding K\"ahler form $\omega$ and the almost complex structure $J_\psi$ with $\psi = d\omega$. By Theorem  \ref{startingpoint} we have that $J = \pm J_\psi = J_{\pm \psi}$ and therefore the NK structure ($g,J$) is uniquely associated to ($f_1,f_2,f_3,f_4,f_5$) or to ($-f_1,-f_2,-f_3,-f_4,-f_5$), where the $f_i's$ are the functions described in the theorem above.
\end{rem}
\bigskip
\section{Cohomogeneity one NK structures that are locally homogeneous}
\setcounter{equation}{0}
\subsection{Local and  global homogeneity}
\begin{prop} \label{homogeneity} Let $(M, g, J)$ be a 6-dimensional strict NK-manifold, admitting  a cohomogeneity one action of   $ \SU_2 \times \SU_2$  preserving the NK structure.  Then $(M, g, J)$ is locally homogeneous if and only if it is locally equivalent  to one of the following compact homogeneous NK spaces: (a) the standard sphere $S^6 = \G_2/\SU_3$;\ (b)
the twistor space $\bC P^3 = \Sp_2/\T^1 \times \Sp_1$;\ (c) the homogeneous space $S^3 \times S^3 = \SU_2^3/(\SU_2)_{\operatorname{diag}}$.
\end{prop}
\begin{pf}  Given $p_o \in   M$, let  $\gg$ be the Lie algebra
of the germs at $p_o$ of  the Killing vector fields of $(M, g)$ preserving $J$,  and denote by   $\gk \subset \gg$ the isotropy subalgebra at $p_o$. By hypothesis
$$\gg \supset \su_2 + \su_2\ ,\qquad \gk \supset \gk \cap (\su_2 + \su_2) = \bR\ .$$ We  recall also  that $\gk$ is reductive and naturally embeddable into $\su_3$.  It  follows that it is isomorphic to $\bR$, $\bR^2$, $\su_2$, $\bR + \su_2$ or $\su_3$.
We denote by $G$  the simply
connected Lie group with  $ Lie(G) = \gg$ and by $K \subset  G$  the connected Lie subgroup
with  $Lie(K) = \gk$. \par
We claim that $K$ is closed in $G$.
Suppose not, so that  $\overline \gk = Lie(\overline K) \supsetneq \gk$.
In this case $\gk$ is
an ideal of $\overline \gk$ and therefore, using the reductiveness of $\gk$, there exists an $\ad(\gk)$-invariant  decomposition  $\overline \gk = \gk + \gp$ with $\gp\neq \{0\}$ and $[\gk,\gp] = 0$. \par
Consider the reductive decomposition  $\gg = \gk + V$ of  $\gg$, with $V
\simeq \bC^3$ and the adjoint action $\ad_\gk|_V\subseteq \su(V)$. From the existence of the space  $\gp \neq\{0\}$ commuting with $\gk$, we see that $\gk$ is isomorphic to $\bR$ or to $\su_2$. Since $K$ is not closed, the Lie algebra $\gk$ is not semisimple (see e.g. \cite{He} p.152), hence $\gk \cong \bR$ and therefore $\gk =
\gk \cap (\su_2 + \su_2)$. Hence $K$ coincides with the isotropy subgroup of $S = \SU_2 \times \SU_2$. In particular,
 $K$ is closed in $S$. Being  semisimple,  $S$ is closed in $G$
and   $K$ is closed in $G$, too, contradicting  $K \subsetneq \overline K$.
\par
\smallskip
Being $K \subset G$ closed, $(M, g)$ is locally isometric to an homogeneous NK manifold $(\wt  M = G/K, \wt g,
\wt J)$ (\cite{Sp, Sp1}).  Looking at the
classification of 6-dimensional homogeneous NK manifolds (\cite{Bu}) and using  \cite{PS} Thm. 1.1, one can directly check that an  homogeneous NK manifold admits a cohomogeneity one action of $\SU_2 \times \SU_2$ if and only if it is one of those listed  in the statement. \end{pf}
\subsection{Locally homogeneous NK structures and their associated  functions  $f_i$'s.}
\hfill\par
\subsubsection{The invariant NK structure of $S^6 = \G_2/\SU_3$}\label{sphere} \label{lochom1}
We   identify the standard 6-dimensional sphere   with the unit sphere $S^6 \subset \Im \bO \simeq \mathbb R^7 $  in  the space of  imaginary octonions,  endowed  with the standard metric $g$ induced by     the Euclidean  product $< x, y> = \frac{1}{2}( x \bar y + y \bar x)$.
We recall that  the algebra of octonions $\bO$  can be  defined as the vector space $\bH^2 = \bH \oplus (\bH \cdot \varepsilon)$,  endowed with the   product  rule (see e.g. \cite{Br})
$$(q_1 + q_2 \varepsilon)\cdot (r_1+ r_2 \varepsilon) \=  (q_1 r_1 - r_2^* q_2) + (r_2 q_1 + q_2 r_1^*) \varepsilon\ .$$
The sphere inherits a natural almost complex structure $J$, which is defined as follows: for $p\in S^6$ and $v\in T_pS^6$, we define $J_p (v)  \= \frac{1}{2}(p \cdot v - v\cdot p)$. It is well known  that the pair $(g, J)$ is a $G_2$-invariant strict NK structure with scalar curvature $s = 30$ (\cite{Gr0}).\par
Identifying unit quaternions with elements of $\SU_2$ in a standard way, we have  the following action of
 $G:=\SU_2\times\SU_2$ on $ \Im \bO$
\beq (q_1,q_2)\cdot (a+b\cdot\varepsilon) = (q_1\cdot a\cdot q_1^*) + (q_2\cdot b\cdot q_1^*) \cdot \varepsilon\ .\eeq
This action provides an embedding of $G\subset \G_2 = Aut(S^6,g,J)$ and the orbit space $S^6/G$ is one-dimensional with principal orbit diffeomorphic to $\SU_2 \times\SU_2/T^1_{\operatorname{diag}}$, and  two singular orbits
$$G\cdot i \simeq \SU_2 \times \SU_2/T^1\times \SU_2 = S^2\ ,\quad G\cdot \varepsilon \simeq \SU_2 \times \SU_2/\SU_2{}_{\operatorname{diag}} =  S^3\ .$$
The curve  $\g_t =  \cos t  \cdot i+ \sin t\cdot \varepsilon$ is a normal geodesic for this action, i.e. it is a geodesic orthogonal to all $G$-orbits. Moreover, a basis for $\su_2+  \su_2$,  with the  same properties  of    \eqref{basissu2}, is given by
$$ U = {\tiny \left(\left(\begin{array}{cc}
\frac{i}{2} & 0  \\
0 & -\frac{i}{2}  \end{array}
\right) , \left(
\begin{array}{cc}
\frac{i}{2} & 0  \\
0 & -\frac{i}{2}  \end{array}
\right)\right)} \ ,\quad A =
{\tiny\left(\left(
\begin{array}{cc}
\frac{i}{2} & 0  \\
0 & -\frac{i}{2}  \end{array}
\right) ,  \left(
\begin{array}{cc}
-\frac{i}{2} & 0  \\
0 & \frac{i}{2}  \end{array}
\right)\right)}\ ,$$
$$E_1 = \tiny\left(\left(
\begin{array}{cc}
0 & \frac{1}{2\sqrt{2}}  \\
-\frac{1}{2\sqrt{2}} & 0  \end{array}
\right) ,  0\right)\ ,\quad V_1 =
\left(\left(
\begin{array}{cc}
0 & \frac{i}{2\sqrt{2}}  \\
\frac{i}{2\sqrt{2}} & 0  \end{array}
\right) , 0\right)\ ,$$
\beq \label{basisexpl} E_2 = {\tiny\left(0, \left(
\begin{array}{cc}
0 & \frac{1}{2\sqrt{2}}  \\
-\frac{1}{2\sqrt{2}} & 0  \end{array}
\right) \right)}\ ,\quad V_2 ={\tiny
\left(0, \left(
\begin{array}{cc}
0 & \frac{i}{2\sqrt{2}}  \\
\frac{i}{2\sqrt{2}} & 0  \end{array}
\right) \right)}\ .\eeq
It is now easy to check that the K\"ahler form $\o$ is given by $\o = f_i \o^i$ , where the $2$-forms $\o_i$, $1 \leq i \leq 5$, are defined as in \eqref{invariantforms} and the functions $f_i$ are as follows
$$f_1= -\sin t \ ,\qquad f_2 = \frac{1}{8}(4 - 9 \sin^2 t)\cdot\cos t\ ,\quad
f_3 = -\frac{1}{8}\sin^2 t\cdot  \cos t\ ,$$
$$f_4 = 0\ ,\quad f_5 = \frac{3}{8}\sin^2 t\cdot\cos t\ . $$
It is immediate to check that, for $t \in ]0,\frac{\pi}{2}[$,  these functions  satisfy   Theorem \ref{theoequations} (i) - (iii) with $\mu = \frac{s}{30} = 1$. \par
\bigskip
\subsubsection{The invariant NK structure of $\bC P^3 = \Sp_2/\T^1 \times \Sp_1$}
\label{lochom2}  It is known that there exists an invariant strict NK structure $(g, J)$ on the twistor space
$\Sp_2/\T^1\times \Sp_1 = \bC P^3$, with scalar curvature $s = 60$,  which can be described as
follows.  We consider the $\ad(\gt_1 + \sp_1)$-invariant decomposition $\sp_2 = (\gt_1 +  \sp_1) + \gp_1 +  \gp_2$
where $\gp_1\cong \bR^2$ and $\gp_2 \cong \bH$. The module $\gp_1 + \gp_2$ identifies with the tangent space
of $\bC P^3$ at the origin $o:= [\T^1 \times \Sp_1]$ and the metric $g$ can be described as the unique  $\Ad(\T^1 \times \Sp_1)$-invariant scalar product on $\gp_1+\gp_2$  with the following properties: $g(\gp_1,\gp_2)=0$, it induces
on $\gp_2 = \bH$ the standard Euclidean product $g(q_1,q_2) = Re(q_1^*\cdot q_2)$
and $g(W,W) = \frac{1}{2}$, where $W = \operatorname{diag}(j,0)\in \gp_2$ (see e.g. \cite{Zi}). Finally,   $J$ is defined as the unique $\Sp_2$-invariant  almost complex structure which corresponds to  the multiplication by $-i$ on $\gp_1$ and by $i$ on $\gp_2$. \par
\medskip
The subgroup $G = \Sp_1\times \Sp_1\subset \Sp_2$ acts on $\bC P^3$ with codimension one principal orbits $G$-equivalent to $\SU_2 \times\SU_2/T^1_{\operatorname{diag}}$. A singular orbit is $G \cdot o$  and the curve $\gamma_t = \exp({\left(
\smallmatrix
0 & t  \\
- t & 0  \endsmallmatrix
\right)})\cdot o$ is a normal geodesic for the action. \par
A basis for $\su_2 + \su_2 \subset \sp_2$ as in \eqref{basissu2} is given by the following matrices
in $\sp_2$
$$  U = \operatorname{diag}\left(\frac{i}{2},\frac{i}{2}\right)\ ,\ \   A =
\operatorname{diag}\left(\frac{i}{2},-\frac{i}{2}\right)\ ,\ \     E_1 = \operatorname{diag}\left(\frac{j}{2\sqrt{2}},0\right) \ ,$$
$$ V_1 = \operatorname{diag}\left(\frac{k}{2\sqrt{2}},0\right)\ ,\ \  E_2 =  \operatorname{diag}\left(0,\frac{j}{2\sqrt{2}}\right)\ ,\ \   V_2 =  \operatorname{diag}\left(0,\frac{k}{2\sqrt{2}}\right)$$
As in the previous example, an easy computation shows that
$$f_1 = \sin t\cdot \cos t \ ,\quad  f_2 = \frac{1}{16}(2\sin^2 t - \cos^2 t )\cdot\cos^2 t \ ,$$
$$f_3 =
\frac{1}{16}(2\cos^2t  - \sin^2 t )\cdot\sin^2 t \ ,\quad f_4 =  0\ ,\quad  f_5 = -\frac{3}{16}\sin^2 t\cdot \cos^2 t\ .  $$
For $t\in ]0,\frac{\pi}{2}[$, the functions $-f_i$'s satisfy Theorem \ref{theoequations} (i) - (iii) with $\mu = \frac{s}{30} = 2$ (see Remark \ref{succinilcolina}). \par
\bigskip
\subsubsection{The invariant NK structure of  $S^3 \times S^3 = \SU_2^3/(\SU_2)_{\operatorname{diag}} $}  \label{lochom3} We recall  that   $\SU_2^3/(\SU_2)_{\operatorname{diag}} $ can be naturally identified with  $S^3 \times S^3$ using  the fact that $L:= \SU_2^3$  acts transitively on $S^3\times S^3 \cong \SU_2\times \SU_2$ by the map
$$(g_1,g_2,g_3)\cdot (x_1, x_2) = (g_1x_1g_3^{-1},g_2x_2g_3^{-1})$$
 with isotropy   $(\SU_2)_{\operatorname{diag}}\subset \SU_2^3$.
We denote by $\gh = \{(X, X, X) \ ,\ X \in \su_2\ \}$ the isotropy subalgebra of $\gl$  and we consider the  $\ad(\gh)$-invariant decomposition
$\gl = \gh +\gm$, in which
$\gm=\{\ (X_1,X_2,X_3)\in \gl\ , \
\sum_{i=1}^3 X_i = 0\ \}$.\par
It is known that there exists an  invariant NK structure $(g, J)$ on $ \SU_2^3/(\SU_2)_{\operatorname{diag}} $, with scalar curvature $s = 60$,  which  is defined as follows. The almost complex structure $J$ is (up to a sign) the unique invariant tensor which acts on $\gm$ as follows
$$J(X_1,X_2,X_3) = \frac{1}{\sqrt{3}}(2X_3+X_1,2X_1+X_2,2X_2+X_3). $$
The metric $g$ is the invariant Riemannian metric which induces on $\gm$ the inner product
 $g:= - \frac{1}{24}
(\B\oplus \B \oplus \B)|_{\gm\times\gm}$, where
 $\B$  is  the Cartan-Killing form of $\su_2$. \par
We consider the subgroup $G\subset L$ given by $G = \{(g,h,g)\in L;\ g,h\in \SU_2\}$. Then $G\cong\SU_2\times \SU_2$ acts on $\SU_2^3/(\SU_2)_{\operatorname{diag}} $ with cohomogeneity one; the following curve is easily checked to be a normal geodesic
\beq\label{geoS^3S^3} \g_t = \exp(t N)\cdot o\ ,\qquad N = \sqrt{6}\left(
\left(\begin{matrix} \frac{i}{2} & 0 \\ 0 & \frac{-i}{2}\end{matrix}\right), 0,  \left(\begin{matrix}- \frac{i}{2} & 0 \\ 0 & \frac{i}{2}\end{matrix}\right)\right)\in \gm.\eeq
If we choose the same basis $U,A,E_1,V_1,E_2,V_2$ of $\su_2 + \su_2$ as in section \S \ref{sphere}, we easily see
that the corresponding functions $f_i$ are given by $$ f_1 = -\frac{\sqrt{2}}{3} \ ,\quad  f_2 = \frac{\sqrt{3}}{36}\sin(2 \sqrt{6} t)\ ,\ \  f_3 = 0 \ , $$
$$  f_4 =  0\ ,\quad  f_5 = -\frac{\sqrt{3}}{36}\sin(\sqrt{6}t)\ ,  $$
which satisfy Theorem \ref{theoequations} (i) - (iii) with $\mu = \frac{s}{30} = 2$ for $t\in ]0,\frac{\pi}{2\sqrt{6}}[$.
\bigskip
\subsection{The functions $f_i$'s as isometric invariants}
\medskip
\begin{prop} \label{isometryclasses} Let $(\o = f_i \o^i, \psi)$ and $(\bar\o = \bar f_i \o^i,\bar\psi)$  be  pairs as in \eqref{3.2} satisfying Theorem \ref{theoequations} (i) - (iii) and associated  with NK structures $(g, J)$ and $(\bar g, \bar J)$, respectively, on $M = ]a,b[ \times G/K$.
Let also $ \fff = \sqrt{f_4^2 + f_5^2}$ and $\bar\fff = \sqrt{\bar f_4{}^2 + \bar f_5{}^2}$.
\par
There exists a local isometry $\varphi: \cU \subset M \longrightarrow M$ between $g$ and $\bar g$  if and only if,  on a suitable subinterval $I \subset ]a,b[$ and modulo compositions with suitable shifts of parameters $t \longmapsto  t + c$, the quadruple $(f_1,f_2,f_3,\fff)= \tau(\bar f_1,\bar f_2,\bar f_3,\bar \fff)$ where $\tau$ belongs to the group of transformations  generated by
\beq \label{merda1} \tau_1(x_1(t),x_2(t),x_3(t),x_4(t)) = (- x_1(-t), x_2(-t), x_3(-t), x_4(-t)),\eeq
\beq\label{merda2}\tau_2(x_1(t),x_2(t),x_3(t),x_4(t)) = (-x_1(t),-x_2(t),-x_3(t),x_4(t)),\eeq
\beq \label{merda3} \tau_3(x_1(t),x_2(t),x_3(t),x_4(t)) = (x_1(t),-x_3(t),-x_2(t),x_4(t)).\eeq \end{prop}
\begin{pf} Let $\varphi:\cU\to \cV$  be an isometry, where $\cU,\cV$ are open subsets of $M$, with $\varphi_* g = \bar g$, hence $\varphi_* J = \pm \bar J$. We denote by $\gs$ and $\bar\gs$ the Lie algebras of Killing vector fields on ($\cU,g$) and ($\cV,g'$) respectively, so that the map $\varphi$ induces a Lie isomorphism $\varphi_*:\gs\to \bar\gs$.\par
\medskip
First of all, we assume that $\gs \cong \bar\gs\cong \gg = \su_2 \oplus\su_2$. \par
In this case,   $\varphi$ maps (locally) any $G$-orbit  into a  $G$-orbit and  therefore $\varphi_*(\xi) = \pm \xi$.
Replacing $\varphi$ by $h\circ\varphi$ for a suitable $h\in G$ and up to a reparameterization $t \mapsto t + c$ , we can suppose that $ \varphi(\g_t) = \g_{\pm t}$.
We can always reduce the case  $\varphi(\g_t) = \g_{ t}$,  by possibly considering the map $\s: (t, gK) \longmapsto (-t, gK)$ and the
isometry $\varphi' = \varphi \circ \s$ between $ \s_* g$ and $\bar g$. Notice that  the  quadruple $(f'_1,f'_2,f'_3,\fff')$ associated with $(\s_*g, \s_* J)$ is obtained
from the quadruple of $(g, J)$ by the transformation
\eqref{merda1}. \par
 If we now think of $\varphi_*$ as an automorphism of the Lie algebra $\gg$, we see that $\varphi_*(\gk) = \gk$ because $\varphi$ preserves $\g$. It then follows that $\varphi_*$ maps the centralizer of $\gk$ in $\gg$ onto itself, hence $\varphi_*(A) = \pm A$.  \par
\noindent{\it Case 1:  $\varphi_*(A) = A$\/}.
If $\gn_i$ ($i=1,2$) denotes the linear span of $\{E_i,V_i\}$ in $\gg$, we see that $\varphi_*$ either preserves or exchanges $\gn_1$ and $\gn_2$.
\begin{lemma} If $\varphi_*(\gn_i) = \gn_i$ for $i=1,2$, then the two quadruples coincide.\end{lemma}
\begin{proof} Put $\Phi_i = \varphi_*|_{\gn_i}\in \End(\gn_i)$ for $i=1,2$. Since $\varphi_*(A)=A$, we have that $\Phi_i$ commutes with $\ad(A)|_{\gn_i}$, hence $\det \Phi_i = 1$. With respect to the basis $B_i =\{E_i,V_i\}$ of
$\gn_i$, we can write the matrix of $\Phi_i$ as  $\left(\smallmatrix \cos(\theta_i)& \sin(\theta_i)\\ -\sin(\theta_i)&\cos(\theta_i)\endsmallmatrix\right)$ for some $\theta_i\in \bR$, $i=1,2$. This implies that
$\varphi^*(\omega_j)|_{\g_t} = \omega_j|_{\g_t}$ for $j=1,2,3$, while
$\varphi^*(\o^4)|_{\g_t} = \cos(\theta_1 - \theta_2) \o^4|_{\g_t}$ and $\varphi^*(\o^5)|_{\g_t} = \sin(\theta_1 - \theta_2) \o^5|_{\g_t}$. Since $\varphi^* (f_j \o^j)|_{\g_t} = f_j  \varphi^*(\o^j)|_{\g_t} = \bar f_j \o^j|_{\g_t}$,
the claim follows.\end{proof}
\begin{lemma} If $\varphi_*(\gn_1) = \gn_2$ and $\varphi_*(\gn_2) = \gn_1$, then
$(f_1,f_2,f_3,\fff)=\! \tau_3(\bar f_1,\bar f_2,\bar f_3,\bar \fff)$. \end{lemma}
\begin{proof} Let $\Phi = \left(\smallmatrix 0& \Phi_1\\ \Phi_2&0\endsmallmatrix\right)$ be the matrix
of the endomorphism $\varphi_*|_{\gn}$ w.r.t. the basis $B$ given by $(E_1,V_1,E_2,V_2)$, where $\Phi_i\in O(2)$ for $i=1,2$. Since $\varphi_*$ commutes with $\ad(A)$, we have that $\det\Phi_i=-1$ for $i=1,2$. This implies that
$\varphi^*\o^1|_{\g_t} = \o^1|_{\g_t}$ and $\varphi^*\o^2|_{\g_t} = -\o^3|_{\g_t}$,
$\varphi^*\o^3|_{\g_t} = -\o^2|_{\g_t}$, while $\varphi^*\o^4|_{\g_t} = (a\o^4+b\o^5)|_{\g_t}$
and $\varphi^*\o^5|_{\g_t} = (-b\o^4+a\o^5)|_{\g_t}$ for some constant $a,b$ with $a^2 + b^2 = 1$. Since
$f_j  \varphi^*(\o^j)|_{\g_t} = \bar f_j \o^j|_{\g_t}$, we get our claim.\end{proof}
\noindent{\it Case 2:  $\varphi_*(A) = -A$\/}. We consider  the automorphism $\psi$ of $G = \SU_2 \times \SU_2$ given by $\psi(g_1,g_2) =
(\bar g_2, \bar g_1)$ and the corresponding  diffeomorphism $\hat \psi$ of $G/K$ and of $M$ given by $\hat \psi(t,gK) = (t,\psi(g)K)$. Note that $\psi_* A= -A$ and that
the NK structure ($\hat\psi_*\bar g, \hat\psi_*\bar J$) is represented by the quadruple $\tau_2(\bar f_1,\bar f_2,\bar f_3,\bar \fff)$. By construction, the isometry $\hat\varphi := \hat\psi\circ\varphi$ satisfies the hypothesis of Case 1 and we get our claim.\par
The sufficiency can be dealt with in a similar way.\par
\medskip
We now discuss the case where  $\gs \supsetneq \gg$. By the results in \cite{PS} and Proposition \ref{homogeneity}, this occurs when the NK structure $(g, J)$ and $(\bar g, \bar J)$ are  both locally homogeneous with  $\gs \simeq \gg_2$, $\sp_2$ or $\su_2^3$.
In all these cases, there exists a local isometry $\psi: \cV \subset M \to \cV \subset M$  of $(M, \bar g)$, which preserves  $\bar J$ and  so that $(\psi \circ \varphi)_*(\gg)  = \gg$.  Indeed, in the first two cases there is only
one subgroup locally isomorphic to $G$ in $\G_2$ or $\Sp_2$ up to conjugation. In the third case, two subgroups
of $L:=\SU_2^3$, isomorphic to $\SU_2^2$ and acting by cohomogeneity one on $Q:= \SU_2^3/\SU_2$ are related by an outer automorphism $\sigma$ of $L$ interchanging two factors and any such $\sigma$ induces an automorphism of the NK structure on $Q$.\par
 Hence, there exists a  local equivalence  $\varphi$ between the NK structures  $(g, J)$ and $ (\bar g,\bar J)$ if and only if there exists a (possibly different)  local equivalence of the two NK structures that,  in addition, maps $\gg$ into itself. From this, the conclusion follows from the first part of the proof.
 \end{pf}
\bigskip
\section{The space of solutions of the differential system and non-homogeneous NK structures}\label{H}
\setcounter{equation}{0}
\medskip
In this section, we study  the space of solutions to the differential problem described  in  Theorem \ref{theoequations} (ii) and (iii). To this purpose,  we  consider the following change of variables.  For any map  $(f_1, f_2, f_3, \fff): ]a,b[ \subset \bR \to \bR^4$ with $f_1 < 0$,  choose  $t_o \in ]a,b[$  and define
$$s(t) \= \int_{t_o}^t \frac{1}{f_1(u)}\ du\ , \qquad g(t) \= \frac{1}{2} \int_{t_o}^t f_1(u) \ du\ .$$
 Being $s' = \frac{1}{f_1} < 0$, we may consider the inverse change of parameter  $t = t(s)$ and set
\beq\label{change} h_1(s) \= g(t(s))\ ,\qquad  h_2(s) \= f_2(t(s)) + f_3(t(s))\ ,$$
$$h_3(s) \= f_2(t(s)) - f_3(t(s))\ ,\qquad  h_4(s) \= 2\ \fff(t(s))\ .\eeq
Using the fact that $\frac{dt}{ds} = f_1(t(s))$,  one gets that  $ h_1'(s) = \frac{1}{2} \left(f_1(t(s))\right)^2$ and the system  \eqref{finalsystem}  turns out to be  equivalent to a  handier differential  system on the $h_i$'s. The results are given in   the following proposition. \par
\begin{prop} \label{spacesolutions}  The change of variables \eqref{change} gives a one-to-one correspondence between  the maps $(f_1, f_2, f_3, \fff): ]a,b[ \to \bR^4$ satisfying (ii) and (iii) of  Theorem \ref{theoequations}  and  the
solutions $(h_1, h_2, h_3, h_4): ]\a,\b[ \to \bR^4$  to  the regular O.D.E. system:
\beq
\label{systemH} \left\{\begin{array}{l}
 h_1''  +  \frac{  2 (h'_1)^2 h_3 + \frac{4}{9\mu} h_4' h_4}
{h_2^2 - h_3^2 - h_4^2} = 0\ ,\cr
\ \cr
h_2'' +  24\mu\ h_1' h_2  = 0\ , \cr
\ \cr
h_3'' -   \frac{  2 (h'_1)^2 h_3 + \frac{4}{9\mu} h_4' h_4}
{h_2^2 - h_3^2 - h_4^2} + 24\mu\ h_1'h_3= 0 \ , \cr
\ \cr
h_4'' + 24\mu\ h_1' h_4 - 4 h_4  = 0\ ,
\end{array}\right.\eeq
 with  initial conditions  $h_i(0)\= a_{i }$,
$h'_i(0) \= b_{i }$  satisfying  the equations
\beq
\label{dataH}a_1 = 0\ , \qquad \cI(a_2, a_3, a_4, b_1, b_2, b_3, b_4 ) = 0_{\bR^4}\ ,\eeq
where  the map $\cI = (\cI^1, \dots, \cI^4): \bR^7 \to \bR^4$ is  defined by
\beq
\label{integrals} \left.\begin{array}{l}
\cI^1 = 12\mu (a_2^2 - a_3^2 - a_4^2) + b_1 + b_3 \ ,\cr
\ \cr
\cI^2 = 4 a_4^2 + b_2^2 - b_3^2 - b_4^2   - b_1^2  - 2 b_3 b_1\ ,\cr
\ \cr
\cI^3 = a_2 b_2 - a_3 b_3 - a_4 b_4 -a_3  b_1 \ , \cr
\ \cr
\cI^4 = \frac{9\mu}{2} b_1 (a_2^2 - a_3^2 - a_4^2) + a_4^2\end{array}\right.\eeq
and so that, for any $t$,  the first derivative $b_j(t) = h'_j(t)$ satisfies the inequalities
\beq \label{ineq}\ b_1 > 0\ ,\qquad b_2^2 - b_3^2   -   b_4^2  - b_1^2 - 2 b_1 b_3< 0\ .\eeq
\end{prop}
\begin{pf} One can directly check that under the change of variables \eqref{change} the system
 \eqref{finalsystem} is equivalent to
\beq \label{A2} h_2'' + 24\mu\ h_1' h_2 = 0\ ,\ \ h_3'' + h_1'' + 24\mu\ h_1'h_3 = 0 \ ,\ \  h_4'' + 24\mu h_1' h_4 - 4 h_4 = 0\ , \eeq
\beq \label{A4}  h_3' + h_1' + 12\mu\ (h_2^2 - h_3^2 - h_4^2) = 0\ .\eeq
By the proof of Theorem  \ref{theoequations}, we also know that  \eqref{spippolo} is satisfied at $t_o$ if and only if it is satisfied for any $t$. So, by the same arguments, we have that the condition  \eqref{spippolo} is equivalent to the differential equation
\beq \label{A5}  (h_2')^2  - (h_3')^2  - (h_4')^2  - (h_1')^2 - 2 h_1' h_3' +  4\ h_4^2= 0\ .\eeq
Differentiating \eqref{A4} and subtracting $\eqref{A2}_2$,  we get
\beq \label{A4'} h_2h_2' - h_3h_3' - h_4h_4' - h_1' h_3 =  \frac{1}{2} (h_2^2 - h_3^2 - h_4^2)' - h_1' h_3 = 0\ . \eeq
Now, we differentiate  \eqref{A4'} and replace the expressions for the second derivatives $h''_i$, $i = 2,3,4$,  determined by  \eqref{A2}, so to obtain the equation
\beq \label{5.7} (h_2')^2 - (h_3')^2 - (h_4')^2 - h_1' h_3'  - 24\mu\ h_1' (h_2^2 - h^2_3 -  h^2_4) -  4 h^2_4   = 0\ .\eeq
Then, subtracting \eqref{A5}  from \eqref{5.7}, we have
\beq (h_1')^2 +  h_1' h_3' - 8\ h_4^2  - 24\mu\ h_1' (h_2^2 - h^2_3 -  h^2_4)  = 0  \eeq
and, using \eqref{A4}, we get
\beq \label{A6} \ \ h_1'\ (h_2^2 - h_3^2 - h_4^2) +  \frac{2}{9\mu}\ h_4^2 = 0\ .\eeq
Finally, differentiating \eqref{A6} and using \eqref{A4'} we obtain the differential equation
\beq \label{A7} h_1''  +  \frac{  2 (h'_1)^2 h_3 + \frac{4}{9\mu} h_4' h_4}
{h_2^2 - h_3^2 - h_4^2} = 0\ ,\eeq
which together with \eqref{A2} gives the regular system \eqref{systemH}.  From the way the equation \eqref{A7} has been derived, it is clear that the equations \eqref{A4}, \eqref{A5}, \eqref{A4'} and \eqref{A6} are identically satisfied if and only if they are true at $s = 0$. They can therefore be considered as algebraic conditions on the initial conditions, to which one has to add the condition $h_1(0) = 0$ due to the definition of $h_1$. This  gives \eqref{dataH}.    Since the inequalities \eqref{ineq} correspond to  Theorem \ref{theoequations} (ii), the proposition follows.
\end{pf}
\begin{rem} \label{remarko} The proof of the previous theorem shows that the functions $\cI_i$, $1 \leq i \leq 4$, are first integrals of the system \eqref{systemH} and that the curve $f = (f_1, f_2, f_3, \fff)$ satisfying (ii) and (iii) of  Theorem \ref{theoequations}  with $\mu = 1$ can be identified with a
curve in the subset
\beq \label{N} N := \cI^{-1}(0) \cap \{\ b_1 > 0\ \}\cap \{\ b_2^2 - b_3^2 - b_4^2< 0\ \} \subset \bR^7\ .\eeq
Using Proposition \ref{isometryclasses} and the definitions of the functions $h_i$'s, one can check that two maximal solutions of  \eqref{systemH} correspond to locally isometric structures $(g, J)$,   $(\overline g, \overline J)$ (with $g$ and $\overline g$ possibly non positively defined) if and only if their images in $N$  are  equal up to a transformation in the group $T$ of transformations of $\bR^7$ generated by the elements
\beq \label{trans} \tau_1(a_2,  \dots,   a_4, b_1,  \dots,  b_4) =  ( -a_2, a_3, a_4,  b_1, -b_2, b_3,  b_4)\ , $$
$$ \tau_2(a_2,  \dots,   a_4, b_1,  \dots,  b_4) =  ( -a_2, -a_3, a_4, - b_1, -b_2, -b_3,  b_4)\  .\eeq
In particular, we have that a solution of \eqref{systemH} corresponds to a locally homogeneous NK structure if and only if its trace in $N$ is contained  in one of the three curves, corresponding to the $f'_i$'s described in \S \ref{lochom1} -  \ref{lochom3}, up to  transformations in $T$.
\end{rem}
\begin{theo} There exists a 2-dimensional family of non-isometric,  non locally homogeneous, strict  NK structures on $M = ]- \varepsilon, \varepsilon[ \times G/K$ for a suitable  $\varepsilon > 0$.
\end{theo}
\begin{pf} Consider the initial data for the differential problem \eqref{systemH}
\beq \label{initialdata} x_o=  (a_1, a_2, a_3, a_4, b_1, b_2, b_3, b_4) = \frac{1}{36} ( 0, \sqrt{3}, \sqrt{3},
\sqrt{6}, 4, 0, 0, - 2 \sqrt{2})\ ,\eeq
whose solution $H(t)$ corresponds to the homogeneous NK structure described  in  \S  \ref{lochom3}. The subset $N$ defined in \eqref{N} can be easily shown to be a 3-dimensional smooth submanifold  in a suitable neighborhood $\cU$ of  $x_o$. We can also suppose that $\cU \cap \t(\cU) = \emptyset$ for any $\tau\in T$. Since $H(t)$ defines a positively defined metric $g$, we can shrink  $\cU$ so that  any point $x \in \cU \setminus \operatorname{Trace}(H) $ gives a solution corresponding to a non-locally homogeneous strict NK structure. By Remark \ref{remarko},  any 2-dimensional submanifold,  transversal to the solutions   in $\cU$,  gives a 2-parameter family of initial data corresponding to non-equivalent NK structures.
\end{pf}
\section{Cohomogeneity one NK manifolds with one singular orbit $S^3$}\setcounter{equation}{0}
Let $G = \SU_2 \times \SU_2$ and consider its natural cohomogeneity one action  on $M \= G\times_H V$, where the subgroup $H$ is the diagonal subgroup $(\SU_2)_{\operatorname{diag}}$ and the vector space $V\cong \bR^3$ is $H$-isomorphic to the Lie algebra $\gh$ acted on by $H$ via the standard adjoint representation. The manifold $M$ is clearly diffeomorphic to $TS^3 \cong S^3 \times \bR^3$ and it can be realized as a tubular neighborhood of the singular orbit $S\cong S^3$ in $S^6$ or $S^3\times S^3$ endowed with the $G$-manifold structures  described in \S  4.2.1 and \S 4.2.3.

In this section we classify (up isometries) all $G$-invariant strict NK structures  defined on $M$.\par
Let $E_\pm, V_\pm\in\gg$ be  defined by
$$ E_\pm \= E_1 \pm E_2,\qquad V_\pm \= V_1 \pm V_2$$
and note that $\gh = \Span\{ U, E_+, V_+\}$, while $\gm \= \Span\{A, E_-, V_-\}$ is the orthogonal complement of $\gh$ in $\gg$. We also fix the curve $\gamma_t \= [(e, tU)]$ in $M$ with $t\in \bR$ and for any $2$-form $\omega$ on $M$ we denote by $\omega_t \= \omega|_{\gamma_t}$ its restriction to $\gamma_t$. We know that for any invariant $2$-form $\omega$ on $M\setminus S$ the restriction $\omega_t$ is of the form $\omega_t = \sum_{1}^5 f_i \omega^i_t$, where the $\omega^i$'s are defined in \eqref{invariantforms}.
\begin{prop}\label{ext1} A $G$-invariant $2$-form $\o$ on $M\setminus S$ corresponding to $\omega_t = \sum_{1}^5 f_i \omega^i_t$ on $\gamma_t$, $t\neq 0$, admits a smooth extension on the whole $M$ if and only if all $f_i's$ extend smoothly at $t = 0$ and the following conditions are satisfied: denoting $\alpha_i \= f_i(0)$, $\beta_i \= f_i'(0)$, $1\leq i\leq 5$
\begin{itemize}
\item[i)] $f_1$, $f_4$ are even and $f_2,f_3,f_5$ are odd; in particular $\a_2=\a_3=\a_5 = 0$;
\item[ii)]
\beq \beta_3 = \frac{1}{2}\alpha_1 + \b_2,\quad \b_5 = - \frac{1}{4} \a_1 - \b_2, \quad \a_4 = 0.\eeq
\end{itemize}
Moreover if $\o$ extends on $M$, we have that $(\o_{\g_0})^3 \neq 0$ if and only if $\a_1 \neq 0$.
\end{prop}
\begin{pf} Let $t,x,y$ be the cartesian coordinates on $V$ determined by the basis $\left(\frac{1}{\sqrt{2}}U,E_+,V_+\right)$ and notice that $\{\frac{\partial}{\partial t},
\frac{\partial}{\partial x}, \frac{\partial}{\partial y}, \wh{A}, \wh{E_-},\wh{V_-}\}$ is a frame field in a neighborhood of $p_o \=\gamma_0$ with dual coframe
$\{dt,dx,dy, A^*, E_-^*, V_-^*\}$.
We have that
$$\wh{E_+}|_{\gamma_t} = - \frac{t}{\sqrt{2}}\left. \frac{\partial}{\partial y}\right|_{\g_t}, \quad \wh{V_+}|_{\gamma_t} =  \frac{t}{\sqrt{2}}\left. \frac{\partial}{\partial x}\right|_{\g_t}$$
and
therefore (with the notations as in \eqref{invariantforms})
$$E_1^* = - \frac{\sqrt{2}}{t} dy + 2 \wh E_-^*, \quad V_1^* =  \frac{\sqrt{2}}{t} dx + 2 \wh V_-^*,$$
$$E_2^* = - \frac{\sqrt{2}}{t} dy - 2 \wh E_-^*, \quad V_2^* =  \frac{\sqrt{2}}{t} dx - 2 \wh V_-^*.$$
Using \eqref{invariantforms},
one can find that
$$\omega_t = f_1\ dt \wedge A^* + f_2\ \left(E_-^* - \frac{\sqrt{2}}{t} dy\right) \wedge\left(V_-^* + \frac{\sqrt{2}}{t} dx\right) +$$
$$+  f_3\ \left(- E_-^* - \frac{\sqrt{2}}{t} dy\right) \wedge \left( - V_-^* + \frac{\sqrt{2}}{t}dx\right) + \frac{2\sqrt{2}}{t}f_4\ \left( - E_-^*\wedge dy +
V_-^*\wedge dx\right) + $$
$$ + f_5\ \left(-2 E_-^*\wedge V_-^* + \frac{4}{t^2} dx\wedge dy\right) = $$
\beq\label{*} = 2\ \frac{f_2 + f_3 + 2f_5}{t^2}\ dx\wedge dy + \ (f_2 + f_3 - 2f_5)\ E_-^* \wedge V_-^*+\eeq
$$ + f_1\ dt \wedge A^* + \sqrt{2}\frac{f_3 - f_2}{t}\ (dx\wedge E_-^* + dy\wedge V_-^*) + \frac{2\sqrt{2}}{t}f_4\ ( -dx\wedge V_-^*
+ dy \wedge E_-^*).$$
The restriction to $V\setminus\{0\}$ of the $G$-invariant form $\omega$ on $M\setminus S$ gives a $H$-equivariant map
$$\tilde \omega : V\setminus\{0\} \longrightarrow \Lambda^2(V^* + \gm^*) \cong V + \gm + V^*\otimes \gm^*$$
and $\wt \omega$  smoothly extends to the whole $V$ if and only if each component  $\wt\o^V$, $\wt\o^\gm$, $\wt\o^{V^*\otimes \gm^*}$ does.\par
We consider the component $\tilde\o^V$. Under suitable identifications, we see that $\wt\o^V$ is a $\SO_3$-equivariant map such that
$\wt\o^V(t,0,0) = (\phi(t),0,0)$, where $\phi(t)\= 2\ \frac{f_2 + f_3 + 2f_5}{t^2}$ for $t\neq 0$. It is then easy to see that $\wt\o^V$ extends smoothly on the whole $V$ if and only if $\phi$ extends to a smooth \emph{odd} function on $\bR$. This means that $f_2 + f_3 + 2f_5$ extends as an odd function with
\beq \label{prima}\b_2 + \b_3 + 2 \b_5 = 0.\eeq
The condition on $\tilde\o^\gm$ is similar and its extendability is equivalent to the extendability of $f_2 + f_3 - 2f_5$ as on odd function. \par
We now identify the $H$-modules $V$ and $\gm$ by means of the map $U\mapsto A$, $E_+\mapsto E_-$, $V_+ \mapsto V_-$ and we further split $V^*\otimes \gm^* = S^2(V^*)\oplus V$ and $\wt\o^{V^*\otimes \gm^*} = \wt\o_1 + \wt\o_2$ accordingly. From \eqref{*} we have that the condition on $\wt\o_2$ is equivalent to say that $\frac{f_4}{t}$ extends
as an odd function, i.e. $f_4$ extends evenly with $f_4(0) = 0$.\par
As for the extendability of $\wt\o_1$, we first write $\wt\o_{1 t} \= \wt \o_1|_{\g_t}$ as
$$\wt \o_{1 t} = \frac{1}{\sqrt{2}}\left(f_1\ dt^2 + \frac{2(f_3 - f_2)}{t}\ (dx^2 + dy^2)\right).$$
From this and the invariance under the symmetry $t\mapsto -t$, we see that both $f_1$ and $\frac{f_3 - f_2}{t}$ must extend smoothly to even functions on the whole $\bR$. This implies that $f_2 - f_3$ extends as an odd smooth function. By previous remarks we deduce that $f_2,f_3$ and $f_5$ extend as odd
smooth functions.\par
We now determine $\wt\o_1$ explicitly at any point $p= (t,x,y)\in V$ with $x^2 + t^2 \neq 0$. If $\theta,\phi$ are
defined as
$$\sin\theta = \frac{y}{\r},\ \  \cos\theta = \frac{\sqrt{t^2+x^2}}{\r}\ \ \text{where}\ \r = |p|\ , \quad\sin\phi = \frac{x}{\sqrt{t^2+x^2}},\ \  \cos\phi = \frac{t}{\sqrt{t^2+x^2}},\quad $$
 then the transformation
$$B \= \left(\begin{matrix} \cos\theta\cdot\cos\phi & \cos\theta\cdot\sin\phi&\sin\theta\\
-\sin\phi& \cos\phi&0\\
 -\sin\theta\cdot\cos\phi & -\sin\theta\cdot\sin\phi&\cos\theta\end{matrix}\right)$$
maps $p$ into  $(\r, 0,0)$.
Hence $\wt\o_1(p) = \wt\o_1(B^{-1}(\r,0,0)) = B^{-1}\cdot\wt\o_1(\r,0,0)$ is equal to
$$\wt\o_1(p) = \frac{1}{\sqrt{2}}\left[ f_1(\r)\ (\cos\theta\cdot\cos\phi\ dt +\cos\theta\cdot\sin\phi\ dx  +\sin\theta\ dy)^2 +\right. $$
$$ + \l(\r)\ \left( (-\sin\phi\ dt + \cos\phi\ dx)^2 +
(-\sin\theta\cdot\cos\phi\ dt -\sin\theta\cdot\sin\phi\ dx + \cos\theta\ dy)^2\right)] = $$
$$= \frac{1}{\sqrt{2}} \left(f_1(\r) \frac{t^2}{\r^2} + \l(\r)\left(1 - \frac{t^2}{\r^2}\right)\right) dt^2 +  \frac{1}{\sqrt{2}}\left(f_1(\r)\ \frac{x^2}{\r^2} + \l(\r)\left (1 - \frac{x^2}{\r^2}\right)\right) dx^2   +  $$
$$ + \frac{1}{\sqrt{2}}\left(f_1(\r)\ \frac{y^2}{\r^2} + \l(\r) \left(1 - \frac{y^2}{\r^2}\right)\right)  dy^2 +$$
$$ + \sqrt{2} \frac{tx}{\r^2}\left(f_1(\r) - \l(\r)\right) dt\odot dx +  \sqrt{2}\frac{ty}{\r^2}\left(f_1(\r) - \l(\r)\right) dt\odot dy + $$
$$ + \sqrt{2}\frac{xy}{\r^2}\left(f_1(\r) - \l(\r)\right)\ dx\odot dy, $$
where $\l = \frac{2(f_3-f_2)}{\r}$. Since both $f_1$ and $\l$ extend as even smooth functions
we see that $\wt\o_1$ extends smoothly if and only if
\beq\label{seconda} \a_1 = f_1(0) = \l(0) = 2(\b_3 - \b_2).\eeq
From this and
\eqref{prima}, condition   (ii) follows.  Finally, from (i) and (ii), one has that $\o_0 = \a_1 \left( dt \wedge A^* + \frac{\sqrt{2}}{2}(dx \wedge E^*_- + dy \wedge V_-^*)\right)$, from which last assertion follows.\end{pf}
We now study the case when the two-form $\o$ is the K\"ahler form of a strict NK structure.
\begin{prop}\label{ext2} Let $\o$ be a $G$-invariant two form on $M\setminus S$ and $\o_t = \sum_i f_i\o^i$ its restriction to $\g_t$. Suppose that $\o$ is the K\"ahler form of a strict NK structure and that it extends smoothly to the whole $M$. If $ f_1(0)\neq 0$, then
\begin{itemize}
\item[i)]   $f_4 = 0$ and $ f'_5(0) \neq 0$;
\item[ii)]  the $3$-form $d\o$ is stable at all points of $M$.
\end{itemize}
\end{prop}
\noindent\emph{Proof.}\ First of all, notice  that  if $\o$ extends smoothly,  then  the $f_i$'s satisfy the system \eqref{finalsystem} for any  $t\in \bR$.  Now,  to prove i),
recall that,  by the proof of Lemma \ref{lemma2},  one of the  ratios  $\frac{f_4}{f_5}$ or $\frac{f_5}{f_4}$  is constant. Since $f_4$ and $f_5$  extend to an even
and an odd function, respectively, it follows that either $f_4 =  0$ or $f_5 =  0$. \par

Assume  that $f_5  =  0$. Then, by Proposition \ref{ext1}, the even function
 $f_4$ satisfies the differential equation \eqref{finalsystem}$_3$ with
$f_4(0)=0$ and $f'_4(0) = 0$. This implies that  $f_4\equiv 0$ contradicting the fact that $f_4^2 + f_5^2 \neq 0$ by Theorem \ref{theoequations} (i). So $f_4 = 0$, $f_5 \neq 0$ and, by the same argument, $f'_5(0) \neq 0$.\par
For (ii), recall that by Proposition  \ref{ext1} the form $(\o_t)^3$ is a volume form
in $T_{\g_t}M$ for every $t$ and that the $G$-invariant 3-form $d\o$ is stable at all points of $M$ if and only if $P'(d\o|_{\g_t}) < 0$ for any $t$, where  $P'$ is the polynomial map on 3-forms, defined in \S  \ref{preliminary1} on the base of  the  volume form $\t'_t = (\o_t)^3$. This volume form can be expressed in terms of   $\t_t = \o^1_t \wedge \o^2_t \wedge \o_t^3$  by
 $$\t'_t = 6 f_1(f_2 f_3 - f_5^2) \t_t\ .$$
Hence,  from \eqref{star} and  the  expression \eqref{pol} for the  polynomial map $P$, determined using the volume form $\t_t$,  for $ t > 0$ we have
$$P'(d\o|_{\g_t}) = - \frac{f_5^2}{9^3\ f_1^2\ (f_2f_3 - f_5^2)^2} \cdot \left((f'_5)^2 - \left(f'_2 + \frac{f_1}{4}\right) \left(f'_3 - \frac{f_1}{4}\right)\right)\ .$$
Now, $P'(d\o|_{\g_t}) < 0$ for $t\neq 0$ by \eqref{firstcondition} and therefore we only need to prove that $P'(d\o|_{\g_0})) < 0$. This follows from the fact that using \eqref{firstcondition} and Prop. \ref{ext1},
\beq\begin{matrix}P'(d\o|_{\g_0})) &=& \lim_{t \to 0^+} P'(d\o|_{\g_t})) = - \lim_{t \to 0^+}\frac{4 f_5^4}{9^3f_1^4 (f_2f_3 - f_5^2)^2} = \\
&=& - \frac{4(f_5'(0))^4}{9^3f_1(0)^4\ (f_2'(0)f_3'(0) - (f_5'(0))^2)^2} = - \frac{4^5\ (f_5'(0))^4}{9^3f_1(0)^8}< 0.\ \quad\qquad\qed\nonumber\end{matrix}\eeq

\begin{cor}\label{cor}
Let $f_i$, $1 \leq i \leq 5$, be smooth functions on some interval $]-a,a[\subseteq \bR$ with
\begin{itemize}
\item[i)] $f_4\equiv 0$ and $f_1$ is even with $f_1 < 0$;
\item[ii)] $f_2,f_3,f_5$ are odd with
$f_3'(0) = \frac{1}{2}f_1(0)  + f'_2(0)$ and  $f'_5(0) = - \frac{1}{4} f_1(0) - f'_2(0)\neq 0$;
\item[iii)] the $f_i's$ satisfy the differential system \eqref{finalsystem} together with the algebraic condition \eqref{spippolo}.
\end{itemize}
Then there exists a tubular neighborhood $\cT_\varepsilon = G \cdot \g|_{[0, \varepsilon[}$ of $S$, $0 < \varepsilon \leq a$,  and a $G$-invariant strict NK structure $(g,J)$ on a $\cT_\varepsilon$, whose K\"ahler form $\o$ is the $G$-invariant 2-form associated with  $\o_t = \sum f_i\o^i$ at the points $\g_t$.
\end{cor}
\begin{pf} The initial conditions $f_1(0) < 0$ and $f_5'(0) \neq 0$ imply that the functions $f_i$'s satisfy  (i) and (ii) of Theorem \ref{theoequations} on a suitable subinterval $]-\varepsilon,\varepsilon[$. Therefore the two-form $\o$, corresponding to $\o_t = \sum f_i\o^i$, $t \neq 0$,  defines a strict NK structure $(g,J)$ (with $g$  possibly not positive definite) on  $\cT_{\varepsilon}\setminus S$. By  Propositions  \ref{ext1} and \ref{ext2}, $\o$ extends smoothly also at  the singular $G$-orbit $S$ and $d\o$ is  stable everywhere. Hence the corresponding NK structure $(g,J)$ extends smoothly on the whole $\cT_\varepsilon$. \par
It remains to show  that, by  possibly choosing a smaller $\varepsilon> 0$,   the metric $g$ is positive definite on  $\cT_{\varepsilon}$.
Using the notations  of  the proof of Proposition \ref{ext1}, the tangent space  at  $p_o = \g_0$ is identifiable with $T_{p_o}\cT = V\oplus \gm$ as $H$-module.
Note that $U = \g_0'\in V$ and that
$$JU = \lim_{t\to 0} J(\g_t') = \lim_{t\to 0} \frac{1}{f_1}\wh A_{\g_t} = \frac{1}{\a_1}A\in \gm$$
 by \eqref{K} and \eqref{spippolo}.  In particular, we have that
$JV\cap \gm\not=\emptyset$ and, being   $V$ and $\gm$ irreducible $H$-modules, it follows that $JV = \gm$. On the other hand, from  \eqref{*}, we see that
$\o_{p_o}(V,V) =   \o_{p_o}(\gm,\gm) = 0$ and hence that
 $g_{p_o}(V,\gm)=0$. \par
 Since $g$ satisfies  \eqref{2.10},  we have that $g_p(U,U) = \lim_{t\to 0}g(\g_t',\g_t') = 1$. From this,   $H$-invariance and $J$-Hermitianity, it follows that
$g_{p_o}$ is positive definite on $V\times V$ and  $\gm\times \gm$ and  hence on the whole tangent space $T_{p_o}\cT_{\varepsilon}$.  A value  $\varepsilon > 0$,  so that $g$ is positive definite on  $\cT_{\varepsilon}$,   can  be now chosen by a simple continuity argument.\end{pf}
 We now prove the main result of this section.
\begin{theo} There exists a one parameter family of non isometric, $G$-invariant strict NK structures on
$TS^3 \cong S^3 \times \bR^3$. \end{theo}
\begin{pf} In order to find $G$-invariant strict NK structures on $TS^3$, we look for  functions $f_i$'s satisfying (i)-(iii) of Corollary \ref{cor}. Since we  need  solutions of \eqref{finalsystem}   and \eqref{spippolo} with $f_1$ nowhere vanishing, we may replace that system of differential equations and algebraic conditions by those obtained through the change of variables considered in \S\ref{H}, whose notations will be kept throughout the following. Setting $\mu = 2$ and using the change of variables \eqref{change} with $t_o = 0$ and $h_4(s) = 2f_5(t(s))$, we see that looking for  functions $f_i$'s as in Corollary \ref{cor} is equivalent to looking for functions $h_i$, $1 \leq i \leq 4$,  on  $]-\varepsilon,\varepsilon[$ so  that:\par
\begin{itemize}
\item[a)] they are odd;
\item [b)] they satisfy  \eqref{systemH} at any $s\neq 0$;
\item[c)]  they satisfy the initial conditions (here $b_i\= h_i'(0)$)
\beq \label{cond1} b_1 > 0\ ,\qquad b_3 =  -b_1,\qquad
b_2 =  -b_4 \neq 0.
\eeq
\end{itemize}
Indeed, the initial conditions in Corollary  \ref{cor} (i),(ii) are equivalent to
\eqref{cond1}, which in turn imply conditions \eqref{integrals}, since  $a_i = h_i(0) = 0$ for any $1 \leq i \leq 4$. Note also that condition  $b_2 \neq 0$ can be replaced by $b_2 < 0$, since the transformation
$h_2\rightarrow -h_2$ (corresponding to the change $(f_2,f_3) \to (-f_3,-f_2)$) gives rise to equivalent NK structures (see Proposition \ref{isometryclasses}, transformation \eqref{merda3}). \par
\smallskip
Being the system \eqref{systemH} singular at  $s=0$,   a necessary condition for the existence of solutions $h_i$'s satisfying (a) - (c) is that
\beq\label{sing}\lim_{s\to 0} \frac{2(h_1')^2 h_3 + \frac{2}{9} h_4 h_4'}{s} = 2 \left( b_1^3 - \frac{1}{9} b_4^2\right) = 0 \ . \eeq
This means  that  \eqref{cond1} can be replaced by the initial conditions
\beq\label{cond2}
b_2 =   -3 b_1\sqrt{b_1}\ , \qquad
b_3 =  -b_1\ , \qquad
b_4 =   3 b_1\sqrt{b_1}\ , \qquad b_1 > 0
\eeq
and that   the whole  problem corresponds  to   {\it looking for smooth odd functions $h_i$'s satisfying   \eqref{systemH} and  \eqref{cond2} on an interval $]-\varepsilon, \varepsilon[$}. \par
For this,  it convenient to re-write  the smooth odd functions $h_i$'s as follows:
$$h_i = s p_i,\ i=1,2,4,\qquad h_3 = s(p_3 - p_1)$$
where the $p_i$'s are   some {\it even} functions. Then,  \eqref{systemH} and \eqref{cond2} are equivalent to the following system  on  even functions  $p_i$'s
\beq
\label{systemp}\left\{\begin{matrix}
p_1'' + \frac{2}{s}\ [p_1' + \frac{18\ p_1p_1'(p_3 - p_1) + p_4p_4'}{9(p_2^2 - (p_3-p_1)^2 - p_4^2)}] +\phantom{aaaaaaaaaaaaaaaaaaaaaaaaaaaa} \\
\phantom{aaaaaaaaaaaaaaaaaa}+ \frac{2}{s^2}\ \frac{9 p_1^2(p_3-p_1) + p_4^2}{9(p_2^2 - (p_3-p_1)^2 - p_4^2)} +
\ \frac{2 p_1'^2\ (p_3 - p_1)}{p_2^2 - (p_3-p_1)^2 - p_4^2}  = 0\ ,\cr
\ \cr
p_2'' +  \frac{2}{s}\ p_2' + 48\ p_1\ p_2 + 48\ s\ p_1'\ p_2 = 0\ , \cr
\ \cr
p_3'' +  \frac{2}{s}\ p_3' + 48\ p_1\ (p_3 - p_1) + 48\ s\ p_1'\ (p_3 - p_1) = 0\ , \cr
\ \cr
  p_4'' + \frac{2}{s}\ p_4' + 4\ p_4\ (12\ p_1 - 1) + 48\ s\ p_1'\ p_4 = 0\ ,
\end{matrix}\right.
\eeq
satisfying the initial conditions (here,  $c_i \= p_i(0)$, while $p_i'(0) = 0$ by evenness)
\beq\label{cond3}
c_2 =  - 3c_1\sqrt{c_1}\ , \qquad
c_3 =  0\ , \qquad
c_4 =   3c_1\sqrt{c_1}\ ,\qquad c_1 > 0\ .
\eeq
If  we now set $\cP:(-\varepsilon,\varepsilon) \to \bR^4$ with $\cP\= (p_1,\ldots,p_4)$ and $\cQ \= \cP'$, we see that  \eqref{systemp} and \eqref{cond3} are equivalent to a system of the form
\beq\label{system3}\left\{\begin{array}{l} \cP' = \cQ\cr
\phantom{a}\\
\cQ' = \frac{1}{s^2} A(\cP) + \frac{1}{s} B(\cP,\cQ) + C(s,\cP,\cQ)
\end{array}\right.\eeq
with initial conditions
\beq\label{cond4}  \cP(0) = (c_1, -3c_1\sqrt{c_1}, 0, 3 c_1\sqrt{c_1}),\qquad
\cQ(0) = 0\ ,\qquad c_1 > 0,
\eeq
where $A$, $B$ and $C$ are smooth $\bR^4$-valued functions defined on suitable open neighborhoods of $\cP(0)\in \bR^4$, $(\cP(0),0)\in \bR^8$ and $(0,\cP(0),0)\in \bR^9$, respectively.\par
\smallskip
We claim that  there exists  a formal power series solution $\wh \cP$ of \eqref{system3} of the form
\beq\label{ps}\wh \cP(s) = \sum_{n=0} \frac{\cP_{2n}}{(2n)!} s^{2n}, \qquad \wh \cQ(s) = \sum_{n=1} \frac{\cP_{2n}}{(2n-1)!} s^{2n-1}.\eeq
In fact,  if $\wh \cP$,  $\wh \cQ$ are as in \eqref{ps},  for suitable $A_{2n}$, $B_{2n+1}$, $C_{2n} \in \bR^4$, one has that
$$A(\wh \cP(s)) = \sum_{n=0} \frac{A_{2n}}{(2n)!} s^{2n},\qquad B(\wh \cP(s),\wh \cQ(s)) = \sum_{n=0}\frac{B_{2n+1}}{(2n+1)!} s^{2n+1}\ ,$$
$$ C(s,\wh \cP(s),\wh \cQ(s)) = \sum_{n=0} \frac{C_{2n}}{(2n)!} s^{2n}$$
and $\wh \cP$ and $\wh \cQ$ are formal solutions of  \eqref{system3}
if and only if,  for any $n\geq 0$,
\beq \label{oneofmany} \cP_{2n+2} = \frac{A_{2n+2}}{(2n+2)(2n+1)} + \frac{B_{2n+1}}{2n+1} + C_{2n},\eeq
Since
$$A_{2n+2} = dA|_{\cP(0)}\cdot \cP_{2n+2} \quad \operatorname{mod}\ (\cP_0,\cP_2,\cP_4,\ldots,\cP_{2n}),$$
$$B_{2n+1} = \frac{d^{2n+1}}{ds^{2n+1}}B(\wh \cP(s),\wh \cQ(s)) = \frac{d^{2n}}{ds^{2n}}\left(\frac{\partial B}{\partial \cP}\cdot \wh \cP'(s) +
\frac{\partial B}{\partial \cQ}\cdot \wh \cP''(s)\right) =$$
$$ = \left.\frac{\partial B}{\partial \cQ}\right|_{(\cP(0),0)}\cdot \cP_{2n+2} \qquad \operatorname{mod}\ (\cP_0,\cP_2,\cP_4,\ldots,\cP_{2n}),$$
equation \eqref{oneofmany} can be re-written in the form
\beq \cP_{2n+2} = \frac{1}{(2n+2)(2n+1)} dA|_{\cP(0)}\cdot \cP_{2n+2} + \left.\frac{1}{2n+1}\frac{\partial B}{\partial \cQ}\right|_{(\cP(0),0)}\cdot \cP_{2n+2} + D_{2n},\eeq
for a fixed function $D_{2n}$ of $\cP_0,\cP_2,\cP_4,\ldots,\cP_{2n}$. It follows that
there exists a formal series solution of the form \eqref{ps} if and only if for every $n$
$$\mathcal L_{2n}\cdot  \cP_{2n+2} = D_{2n},$$
where $\mathcal L_{2n}$ is the $4\times 4$-matrix
$$\mathcal L_{2n} = Id - \frac{1}{(2n+2)(2n+1)}dA|_{\cP(0)} -\left. \frac{1}{2n+1}\frac{\partial B}{\partial \cQ}\right|_{(\cP(0),0)}.$$
Since
$$dA|_{\cP(0)} = \left(\begin{matrix} 6 & 0& -2&  -\frac{4}{3\sqrt{c_1}}\\ 0 & 0&0&0\\ 0 & 0&0&0\\
0 & 0&0&0 \end{matrix}\right),\qquad
\left.\frac{\partial B}{\partial \cQ}\right|_{(\cP(0),0)} = \left(\begin{matrix} 6&0&0& -\frac{2}{3\sqrt{c_1}}             \\ 0 & 2 & 0 &0 \\
0 & 0 & 2 & 0\\ 0 & 0 & 0 &2\end{matrix}\right)\ ,$$
we have that   that $\det(\mathcal L_{2n}) = \frac{2n^2 - 3n - 8}{(2n+1)(n+1)} \left(\frac{2n-1}{2n+1}\right)^3 \neq 0$ for every $n\geq 0$ and the sequence  $\cP_{2n+2} = \cL_{2n}^{-1} \cdot D_{2n}$ uniquely  determines a formal solution of the form \eqref{ps}.\par
By Malgrange Theorem (\cite{Ma}, Thm. 7.1), the existence of the formal solution $(\wh \cP, \wh \cQ)$ implies the existence of  a smooth solution $(\cP, \cQ)$ of \eqref{systemp}, whose   Taylor expansion at $0$ coincides with $(\wh \cP, \wh \cQ)$. This  solution can be determined  with $\cP$ even: In fact,  given an arbitrary  smooth solution  $\wt \cP$  of \eqref{systemp} with Taylor expansion $\wh \cP$ at $0$,  the even  function
  $\cP(t) = \wt \cP(|t|)$  is  smooth also at $0$ (because $\wh \cP$  has no odd degree monomials) and $(\cP,  \cQ = \cP')$ is a  solution of  \eqref{system3}, as it is  immediately checked. \par
A smooth even solution $\cP$ gives rise to a strict NK structure $(g, J)$ on $TS^3$. Recall that  $g$ is Einstein and hence  real analytic by \cite{DK}. Since the Killing vector field $\wh A$ is  non-vanishing  along  the geodesic $\g_t$ and
  $f_1(t)$ coincides (up to a multiple)  with the norm $||\wh A||_{\g_t}$ (see \eqref{K} and
\eqref{f1}), we have that $f_1(t)$ is real analytic.  By a direct inspection of the system \eqref{systemH}, we also have   that $h_2$, $h_3$, $h_4$ and the corresponding $\cP$ are real-analytic. From this we get  that {\it any  NK structure on $T S^3$ determined by the solutions $\cP$ are uniquely determined by  the initial data} and in particular, by the choice  the parameter $b_1 = c_1 > 0$   in \eqref{cond3}. \par
\smallskip
It remains to show that two NK structures $(g, J)$, $(\overline g, \overline J)$,  corresponding   to initial data    $b_1$, $\overline b_1 >0$, are isometric if and only if $b_1 = \overline b_1$.\par
First of all, we claim that {\it $b_1$ corresponds to a locally homogenous  structure  if and only if $b_1 = \frac{1}{9}$ or $b_1 = \frac{1}{4}$.}
In fact,   $(g, J)$  is  locally homogeneous if and only if there is an isometric equivalence between an open neighborhood $\cU$ of $p_o = \g_0$ and an open subset $\cU'$ of  one of the   compact homogeneous NK spaces $N$, described in Proposition  \ref{homogeneity}. If  $\gg$ is properly identified with a subalgebra of  $\frak{aut}(N)$,  the  isometric equivalence  can be assumed to be  $\gg$-equivariant.  Notice that  the only spaces $N$, admitting a  cohomogeneity one   action of $G$ with a singular orbit  $G/H \cong S^3$,  are $S^6$ and $S^3\times S^3$.
Now, the  claim can be checked using the explicit descriptions of homogeneous NK structures in \S 4.2.  Indeed,  when  $N = S^3\times S^3$ the expressions in \S \ref{lochom3}   immediately  give that    $b_1 = \frac{1}{9}$, while in the case  $N= S^6$ more work is needed. In fact,  one has  to: a) rescale the   metric  considered in  \S \ref{lochom1} in order to have  $\mu = \frac{s}{30} = 2$; b) replace  the
geodesic  $\g$ with $\wt \g(t) = \g(-\sqrt{2}t+\frac{\pi}{2})$, so that $||\wt  \g'|| = 1$ and $\wt \g(0)\in S$; c)  change $J$ into $-J$ so that the new
functions $\wt f_i$  meet all our assumptions on signs. They are  $\wt f_1 = -\frac{1}{\sqrt{2}}\cos(\sqrt{2}t)$ and $\wt f_i(t) = - \frac{1}{2}f_i(-\sqrt{2}t + \frac{\pi}{2})$, $i\geq 2$, and give  $b_1 = \frac{1}{4}$.\par
\smallskip
Consider now initial data $b_1$, $\overline b_1 >0$, corresponding to locally isometric NK structures $(g, J)$, $(\overline g, \overline J)$. By the previous claim,  we may assume that  they  are not
locally homogeneous.  In this case,   if $\phi:\cU \subset TS^3 \to \cU' \subset TS^3 $ is a local  isometry between $g$ and $g'$ on a neighborhood $\cU$ of $p_o = \g_0 \in S \simeq S^3$, we have that    $\phi$ locally maps  $G$-orbits into $G$-orbits.  In particular,  $p_o$ is mapped into a point of
the singular $G$-orbit $S$ and, by  composition with  some element of $G$,   we may assume that $\phi(p_o) = p_o$ and $\phi(\g_t) = \g_t$
for every $\g_t\in \cU$. From this and  Proposition  \ref{isometryclasses}, it follows that the quadruples $(f_1,f_2, f_3,\fff)$ and
$(\overline f_1, \overline f_2,\overline f_3, \overline \fff)$,  corresponding to $(g, J)$,  $(\overline g, \overline J)$,    coincide. This implies that $b_1 = \overline b_1$, and the proof is concluded.\end{pf}

\bigskip\bigskip
\font\smallsmc = cmcsc8
\font\smalltt = cmtt8
\font\smallit = cmti8
\hbox{\parindent=0pt\parskip=0pt
\vbox{\baselineskip 9.5 pt \hsize=3.1truein
\obeylines
{\smallsmc
Fabio Podest\`a
Dip. di Matematica "U.Dini"
Universit\`a di Firenze
Viale Morgagni 67/A
I-50134 Firenze
ITALY}

\medskip
{\smallit E-mail}\/: {\smalltt podesta@math.unifi.it
}
}
\hskip 0.0truecm
\vbox{\baselineskip 9.5 pt \hsize=3.7truein
\obeylines
{\smallsmc
Andrea Spiro
Scuola di Scienze e Tecnologie
Universit\`a di Camerino
Via Madonna delle Carceri
I-62032 Camerino (Macerata)
ITALY
}\medskip
{\smallit E-mail}\/: {\smalltt andrea.spiro@unicam.it}
}
}

\end{document}